\documentclass[11pt]{article}
\usepackage{amssymb}
\usepackage{mathrsfs}
\usepackage{latexsym,amsmath,amssymb,amsfonts,epsfig,graphicx,cite,psfrag}
\usepackage{eepic,color,colordvi,amscd}
\usepackage{ebezier}
\usepackage{verbatim}
\usepackage{subfigure}
\usepackage{enumitem}
\usepackage{algorithm}
\usepackage{algorithmic}
\usepackage{hyperref}
\usepackage{amsthm}
\usepackage{longtable}
\usepackage{comment}
\usepackage{color}
\usepackage{graphicx} 

\usepackage{listings}
\usepackage{xcolor}

\definecolor{mygreen}{rgb}{0,0.6,0}

\lstdefinestyle{mystyle}{
    language=Matlab,
    basicstyle=\ttfamily\small,
    numbers=left,
    numberstyle=\tiny\color{gray},
    stepnumber=1,
    numbersep=5pt,
    backgroundcolor=\color{white},
    showspaces=false,
    showstringspaces=false,
    showtabs=false,
    frame=single,
    tabsize=4,
    captionpos=b,
    breaklines=true,
    breakatwhitespace=false,
    commentstyle=\color{mygreen},
    keywordstyle=\color{blue},
    morecomment=[l]{\%},
    morekeywords={mymodel,add_affine_constraint,sparse},
}
\usepackage{epstopdf}
\definecolor{blue}{rgb}{0,0,0.9}
\definecolor{red}{rgb}{0.9,0,0}
\definecolor{green}{rgb}{0,0.9,0}

\theoremstyle{plain}
\newtheorem{remark}{Remark}

\newtheorem{assumption}{Assumption}

\newtheorem{theorem}{Theorem}

\usepackage{longtable}
\def\qed{\hfill \rule{4pt}{7pt}}

\def\H{\mathcal{H}}
\def\<{\big\langle}
\def\>{\big\rangle}
\def\E{\mathcal{E}}
\def\M{\mathcal{M}}

\def\A{\mathcal{A}}
\def\B{\mathcal{B}}

\def\R{\mathbb{R}}
\def\X{\mathbb{X}}
\def\Y{\mathbb{Y}}

\def\Re{{\rm Rtr}}

\def\N{\mathbb{N}}

\def\S{\mathbb{S}}

\def\dist{{\rm dist}}
\def\dd{{\rm diag}}
\def\P{\mathcal{P}}

\def\Q{\mathcal{Q}}

\def\({\left(}
\def\){\right)}

\let\svthefootnote\thefootnote
\newcommand\blankfootnote[1]{%
	\let\thefootnote\relax\footnotetext{#1}%
	\let\thefootnote\svthefootnote%
}
\def\V{\mathbb{V}}
\def\cV{\mathcal{V}}

\usepackage{amsopn}

\begin{document}

	\title{A preconditioned augmented Lagrangian method for solving semidefinite programming problems} 
		
	\author{Tianyun Tang
\thanks{Institute of 
Operations Research and Analytics, National
         University of Singapore,
         Singapore
         117602 ({\tt ttang@u.nus.edu}).
         }, \quad 
	 Kim-Chuan Toh\thanks{Department of Mathematics, and Institute of 
Operations Research and Analytics, National
         University of Singapore, 
       Singapore
         119076 ({\tt mattohkc@nus.edu.sg}).  The research of this author is supported
by the Ministry of Education, Singapore, under its Academic Research Fund Tier 3 grant call (MOE-2019-T3-1-010).}
	 }
	
	\date{\today}
	
\maketitle


\begin{abstract}
In this work, we propose a preconditioned augmented Lagrangian method (ALM) for solving semidefinite programming (SDP) problems. The preconditioner is implemented via a weighted penalty function in the ALM subproblem, with the weight matrix derived from the projection operator onto the tangent space of the feasible region. This simple yet effective modification significantly accelerates ALM, particularly for ill-conditioned SDPs. By combining the preconditioned ALM with our previously developed feasible method, SDPF [T.~Tang and K.-C.~Toh, SIAM J. Optim., 34 (2024), pp.~2169–2200], we develop SDPF+, an SDP solver capable of handling convex problems with possibly nonlinear objective functions. Extensive numerical experiments demonstrate the efficiency and robustness of SDPF+, showing that it can generally outperform other solvers on large-scale SDPs whose optimal solutions exhibit low-rank structure.
\end{abstract}

\bigskip
\noindent{\bf keywords:} Semidefinte programming, Augmented Lagrangian method, Preconditioner
\\[5pt]


\section{Introduction}\label{sec:Intro}

Let $\S^n$ be the set of $n\times n$ real symmetric matrices, $\S^n_+$ the set of $n\times n$ symmetric positive semidefinite matrices, and $[m]:=\{1,2,\ldots,m\}.$ Consider the following standard form linearly constrained SDP problem: 
\begin{equation}\label{SDP}
\min\left\{ f(X):\ \A(X)=b,\ X\in \S^n_+ \right\} \tag{SDP} 
\end{equation}
where $f:\S^n\rightarrow \R$ is a convex continuously differentiable function, $\A:\S^n\rightarrow \R^m$ is a linear map: $\A(X):=[\<A_1,X\>;\ldots;\<A_m,X\>]$ such that $A_i\in \S^n$ for any $i\in [m]$ and $b\in \R^m.$ 
For simplicity of presentation, we consider the standard form SDP expressed in
the {form} \eqref{SDP} first. In Section~\ref{Sec:gene}, we will consider more general forms of SDP. 

For large-scale SDP problems with low-rank optimal solutions, the Burer--Monteiro method~\cite{BM1} reduces the dimensionality by solving~(\ref{SDP}) in the following factorized form:
\begin{equation}\label{SDPR}
\min\left\{ f(RR^\top) \;:\; \A(RR^\top) = b,\; R \in \mathbb{R}^{n\times r} \right\}, \tag{SDPR}
\end{equation}
{with $r\leq n.$} Numerous algorithms have been proposed for solving low-rank SDPs, including the augmented Lagrangian method (ALM)~\cite{BM1}, the alternating direction method of multipliers (ADMM)~\cite{han2024low}, and feasible methods~\cite{tang2024feasible}. Among these, ALM is the most widely used due to its simplicity and strong empirical performance. Several SDP solvers based on ALM have been developed in the literature, such as SDPLR~\cite{BM3}, ManiSDP~\cite{wang2025solving}, HALLaR~\cite{monteiro2024low}, and RiNNAL~\cite{hou2025low}.

ALM is a primal--dual method in which the primal variable is updated by solving the following subproblem:
\begin{equation}\label{SDPLRsub}
\min\left\{ f(RR^\top) - \langle \lambda,\, \A(RR^\top) - b \rangle + \frac{\beta}{2} \|\A(RR^\top) - b\|^2 :\ R \in \mathbb{R}^{n\times r} \right\}, 
\end{equation}
for some given penalty parameter $\beta > 0$ and Lagrange multiplier $\lambda \in \mathbb{R}^m$.  
After obtaining $R^+$ from~\eqref{SDPLRsub}, the multiplier is updated via  
$\lambda^+ \leftarrow \lambda - \beta \big( \A(R^+(R^+)^\top) - b \big)$.  
Since~\eqref{SDPLRsub} is an unconstrained problem, it can be solved by simple first-order methods, potentially with a rank-adaptive strategy~\cite{tang2024feasible}, to obtain a global minimizer. However, when~\eqref{SDPR} is nearly primal degenerate, that is, when the Jacobian of the constraint map is ill-conditioned, solving~\eqref{SDPLRsub} becomes difficult. In this case, the penalty term in~\eqref{SDPLRsub} inherits the ill-conditioning difficulty, causing first order methods such as gradient descent or limited memory BFGS to converge very slowly. In our implementation, the ALM subproblem often requires tens of thousands of gradient descent iterations before the gradient norm becomes sufficiently small to guarantee convergence of the outer ALM iteration. Some ALM-based solvers, such as ManiSDP, limit the maximum number of iterations for the subproblem and update the multiplier once this limit is reached.  
Nevertheless, if the subproblem is not solved with sufficient accuracy, ALM may require a large number of outer iterations or even fail to converge.

Table~\ref{test-illc} reports the performance of several algorithms on the Lov\'asz Theta SDP problem for the graph \texttt{G51} from the Gset dataset\footnote{Dataset available at \href{https://web.stanford.edu/~yyye/yyye/Gset/}{https://web.stanford.edu/\textasciitilde yyye/yyye/Gset/}.}. Although this problem admits a low-rank optimal solution, the interior-point method (IPM) based solver SDPT3~\cite{TTT} significantly outperforms both ManiSDP and SDPLR in terms of computational speed. This advantage stems from the primal degeneracy of the SDP, which severely limits the efficiency of low-rank ALM approaches. Note that the ill-conditioning issue affects not only ALM but also all of the aforementioned algorithms. While IPM is generally more robust for solving ill-conditioned SDP problems, its high computational complexity renders it unsuitable for large-scale instances.

\begin{center}
\begin{small}
\begin{longtable}{|c|ccccc|c|}
\caption{Test results on Lov\'asz Theta SDP problem for graph G51. The condition number of the constraint Jacobian at the optimal solution is around $10^{7}$.}
\label{test-illc}
\\
\hline
problem & algorithm & pfeas& dfeas & comp & fval & ttime \\ \hline
\endhead
G51  & SDPT3 & 3.59e-13 & 0.00e+00 & 8.05e-08 & -3.4899993e+02& 3.08e+01 \\
n=1000& ManiSDP & 7.51e-06 & 2.62e-06 & 2.21e-17 & -3.4902483e+02& 3.60e+03 \\
m=5910 & SDPLR & 9.99e-07 & 4.51e-06 & 1.82e-08 & -3.4900065e+02& 2.74e+03 \\
rank=109 & SDPF & 3.41e-07 & 8.77e-07 & 6.43e-08 & -3.4900017e+02& 5.13e+01 \\
& SDPF+ & 4.67e-07 & 9.59e-07 & 2.36e-09 & -3.4899975e+02& 2.30e+01 \\
\hline
\end{longtable}
\end{small}
\end{center}

\vspace{-11mm}
In order to address the impact of ill-conditioning on ALM, we develop a preconditioned ALM (PALM) for solving~\eqref{SDP}.  
In the ALM subproblem, the constraints $\A(X) = b$ are penalized using a weighted term in the penalty function. The weight matrix, derived from the preconditioner for the projection onto the tangent space of the feasible region, effectively mitigates primal ill-conditioning.  
The resulting algorithm, implemented in the solver SDPF+, substantially accelerates ALM and outperforms SDPT3 on ill-conditioned instances (see Table~\ref{test-illc}). Note that SDPF is a feasible method \cite{tang2024feasible} and employs preconditioning. Table~\ref{test-illc} shows that PALM outperforms the feasible method on this instance.
We further discuss extensions of the PALM to more general SDP problems with inequality and other conic constraints.

The remainder of the paper is organized as follows.  
Section~\ref{Sec-alg} introduces the PALM and establishes its convergence analysis.  
Section~\ref{Sec-solver} presents the algorithmic framework of SDPF+, which employs a two-stage approach that combines the feasible method in \cite{tang2024feasible} with the PALM; details on its usage and data structure are provided in Appendix~\ref{App-solver}.  
Section~\ref{Sec-numer} reports numerical results comparing SDPF+ with other state-of-the-art SDP solvers, with additional tables given in Appendix~\ref{App-table}.  
Finally, Section~\ref{Sec-conc} concludes the paper with a brief summary.

\section{Algorithmic design}\label{Sec-alg}

In this section, we present the PALM. Before introducing the algorithm, we first review the feasible method SDPF, which motivates the design of our preconditioner.

\subsection{Feasible method}

In \cite{tang2024feasible}, we have designed a rank-adaptive feasible method called SDPF to solve \eqref{SDP} via a low-rank factorization approach by solving \eqref{SDPR}. To make the presentation self-contained, we briefly restate the algorithm. It is based on Riemannian optimization \cite{manibook,boumal2023introduction}. Under the linear independence constraint qualification (LICQ)~\cite{nocedal1999numerical}, namely, when the Jacobian of the constraint map has full row rank, the feasible set of \eqref{SDPR} forms a Riemannian manifold.
\begin{equation}
\M := \left\{ R \in \R^{n \times r} : \A(RR^\top) = b \right\},
\end{equation}
whose tangent space at $R$ is given by
\begin{equation}
{\rm T}_R \M := \left\{ H \in \R^{n \times r} : \A(RH^\top + HR^\top) = 0_{m \times 1} \right\}.
\end{equation}
Let ${\rm Proj}_{{\rm T}_R \M} : \R^{n \times r} \to {\rm T}_R \M$ denote the orthogonal projection onto the tangent space, and let $\Re_R : {\rm T}_R \M \to \M$ be a retraction satisfying
\begin{equation}
\Re_R(tH) = R + tH + o(t),
\end{equation}
where $t\in \R$ is a parameter that controls the stepsize. In SDPF, the retraction is implemented via a Gauss--Newton procedure initialized at $R + tH$.
At each iteration, SDPF performs a Riemannian gradient descent step:
\begin{equation}
\label{eq-RGD}
R^+ \leftarrow \Re_R \left( R - t \, {\rm Proj}_{{\rm T}_R \M} \big[ 2 \nabla f(RR^\top) R \big] \right),
\end{equation}
where $t > 0$ is a stepsize parameter.
In each of the update step \eqref{eq-RGD},
the main computational cost lies in evaluating the projection and the retraction, both of which require solving $m\times m$ symmetric positive definite linear systems of the form:
\begin{equation}\label{linsolve}
\frac{1}{2}\A\( \A^*(\lambda)RR^{\top}+RR^{\top}\A^*(\lambda) \)=d,
\end{equation}
where $d\in \R^m$ is a given right hand side vector, and $\lambda\in \R^m$  is the variable. The matrix form of the linear operator in \eqref{linsolve} is given by:
\begin{equation}\label{linsolveM}
M_R:=\begin{bmatrix}
\<A_1R,A_1R\>&\cdots&\<A_1R,A_mR\>\\
\vdots&\ddots&\vdots\\
\<A_mR,A_1R\>&\ldots&\<A_mR,A_mR\>
\end{bmatrix}.
\end{equation}
The positive definiteness of $M_R$ is equivalent to the LICQ condition. Although LICQ may not  hold in general, it can be ensured with probability one by adding a small random perturbation to the vector $b$ in \eqref{SDP} \cite[Theorem~4]{tang2024feasible}.

The primal degeneracy of~\eqref{SDP} corresponds to the ill-conditioning of the matrix $M_R$.  
To efficiently solve the linear systems~\eqref{linsolve} for degenerate SDP problems, we proposed in~\cite{tang2024feasible} an adaptive preconditioning technique that uses  
$M_{R_k}^{-1} = L_{R_k}^{-\top} L_{R_k}^{-1}$,  
where $L_{R_k}$ is the Cholesky factor of $M_{R_k}$ at a selected iteration point $R_k$.  
The preconditioner is updated only when the conjugate gradient (CG) method fails to converge within a small number of iterations. We should mention that the sparsity of the coefficient matrices $\{A_i : i\in[m]\}$ can be exploited when forming $M_R$ and computing its Cholesky factor through the following approaches:
\begin{itemize}
\item For any $i, j \in [m]$, if $A_i$ and $A_j$ do not share common nonzero rows, then $\langle A_i R, A_j R \rangle = 0$. In many SDP problems, most pairs $A_i$, $A_j$ do not have overlapping nonzero rows, thus making $M_R$ very sparse. In such cases, computing the sparse Cholesky decomposition of $M_R$ can often be done with modest cost by using the {\sc Matlab} build-in function \texttt{chol}. 
\item The matrix $M_R$ can be computed efficiently when the $A_i$ are highly sparse.  
If each $A_i$ has at most $K$ nonzero entries, we can first compute $U = A_i R$ at the cost of $O(Kr)$ operations, and then obtain the $i$-th column of $M_R$ as $(M_R)_{:,i} = \A(UR^\top)$ at the cost of $O(mKr)$. The total cost is therefore $O(m^2Kr)$, which is typically acceptable.
\end{itemize}

Although the efficiency of SDPF has been demonstrated through numerical experiments on various SDP problems, its applicability remains limited to instances with a moderate number of constraints (typically $m < 50{,}000$), for the following reasons:
\begin{itemize}
\item Cholesky decomposition may introduce many fill-ins in $L_R$, leading to high memory usage and slow factorization when $m$ is large and $M_R$ lacks a 
conducive sparsity pattern for fast factorization. 

\item Certain SDP problems, such as those arising from polynomial optimization~\cite{lasserre2001global}, involve an extremely large number of constraints (e.g., $m = \Omega(n^2)$). In such cases, computing and storing $M_R$ may exceed available memory.
\end{itemize}

Motivated by the feasible method SDPF,
in the next subsection, we introduce the PALM, which can partially alleviate the above limitations of SDPF.

\subsection{Preconditioned ALM}

The algorithmic framework of our PALM is stated in Algorithm~\ref{alg:PALM}. Note that in the implementation of Algorithm~\ref{alg:PALM}, we adaptively update the rank parameter $r$ to ensure global optimality. More specifically, global optimality of the subproblem is certified by dual feasibility, i.e., the positive semidefiniteness of the dual slack matrix. If this condition fails, we increase the rank of $R$ by an appropriate integer $k$ by appending columns corresponding to the eigenvectors of the $k$ most negative eigenvalues of the dual slack matrix. This provides a descent direction and allows the algorithm to escape the strict saddle point in a higher-dimensional space. See \cite[Section~3.3]{tang2024feasible} for details.

\begin{algorithm}
	\renewcommand{\algorithmicrequire}{\textbf{Input:}}
	\renewcommand{\algorithmicensure}{\textbf{Output:}}
	\caption{Preconditioned ALM for (\ref{SDP})}
	\label{alg:PALM}
	\begin{algorithmic}
		\STATE {\bf Initialization}: Choose $R^1\in \R^{n\times r},\lambda^1\in \R^m,\ \beta>0.$  
		\FOR{$k = 1, 2, \dots$}
		\STATE 1. Choose a positive definite weight matrix $W^k\in \S^m$ ($W^k\approx M_{R^k}^{-1}$ ideally).
		\STATE 2. Solve the following ALM subproblem to get $R^{k+1}$
        \begin{equation}
        \min\left\{ f(RR^\top)-\<\lambda^k,\A(RR^\top)-b\>+\frac{\beta}{2}\|\A(RR^\top)-b\|^2_{W^k}:\ R\in \R^{n\times r}\right\}.
        \label{ALM: step 2}
        \end{equation}
		\STATE 3. $\lambda^{k+1}:=\lambda^k-\beta W^k \( \A(R^{k+1}(R^{k+1})^\top)-b \).$
		\ENDFOR
	\end{algorithmic}
\end{algorithm}

For the subproblem \eqref{ALM: step 2} and the later subproblem \eqref{SDPGALMsub}, we apply gradient descent with a rank-adaptive strategy, which can be viewed as SDPF specialized to unconstrained optimization. The key distinction from solving \eqref{SDPR} directly via SDPF is that no projection or retraction is required, since both subproblems \eqref{ALM: step 2} and  \eqref{SDPGALMsub} are unconstrained. 

Compared to the traditional low-rank ALM, the only modification in Algorithm~\ref{alg:PALM} is the use of a weighted penalty function in the subproblem as a preconditioner, which also influences the update of the Lagrange multiplier in Step~3. This is because the introduction of the weighted penalty term modifies the stationarity condition of the ALM subproblem, and consequently alters the dual ascent step through the incorporation of the weight matrix. Leveraging the equivalence between the primal ALM and the dual proximal point algorithm (PPA)~\cite{rockafellar1976augmented}, ALM with a weighted penalty function can be interpreted as a PPA with a weighted proximal term. While weighted penalty functions and proximal terms have been employed to accelerate ADMM-type splitting methods~\cite{kontogiorgis1998variable, giselsson2014diagonal, bredies2017proximal, pock2011diagonal, liu2021acceleration,sun2025accelerating}, their practical usage in the context of ALM remains relatively unexplored.

We now discuss how to choose the weight in the penalty function to mitigate the ill-conditioning of problem~\eqref{SDPR}. The weight matrix has previously appeared in the dual refinement step of SDPF~\cite[Section~5]{tang2024feasible}; here, we incorporate this idea into ALM and examine it in greater detail. Suppose we fix the preconditioning weight matrix $W^k = W$ in Algorithm~\ref{alg:PALM}. Then PALM can be interpreted as a standard ALM applied to an equivalent reformulation of~\eqref{SDPR}:
\begin{equation}\label{SDPRT}
\min\left\{ f(RR^\top):\ V\A(RR^\top)=Vb,\ R\in \R^{n\times r} \right\}, 
\end{equation}
where $V\in \R^{n\times n}$ satisfies $V^\top V=W.$ Note that (\ref{SDPRT}) is equivalent to (\ref{SDPR}) because $W$ is positive definite and $V$ is invertible. When we use SDPF to solve (\ref{SDPRT}), the linear system in (\ref{linsolve}) will become
\begin{equation}\label{linsolve1}
\frac{1}{2}V\A\( (\A^*(V^\top\mu)RR^{\top}+RR^{\top}\A^*(V^\top\mu)) \)=d_1,
\end{equation}
With the corresponding matrix form being $VM_RV^\top.$ To improve the conditioning of the linear system, an ideal choice is to set $V = L_R^{-1}$, where $M_R = L_R L_R^\top$ is the Cholesky decomposition of $M_R$. This choice ensures that $V M_R V^\top = I_m$, resulting in a well-conditioned system. Consequently, the corresponding weight matrix is given by
\begin{equation}\label{weightW}
W = V^\top V = L_R^{-\top} L_R^{-1} = M_R^{-1},
\end{equation}
which coincides with the preconditioner used in the adaptive preconditioning technique for solving \eqref{linsolve}. This observation suggests that $W^k$ should be chosen to approximate $M_R^{-1}$.

However, this is not directly achievable because $R$ varies during the iterations of the ALM subproblem, and hence $M_R^{-1}$ is not a constant matrix. A reasonable assumption, particularly when the algorithm is close to convergence, is that the primal iterates $R$ remain close to each other within a single ALM subproblem. Under this assumption, $M_{R^k}^{-1}$ -- where $R^k$ is the initial point of the $k$-th ALM subproblem -- can serve as a good approximation of $M_R^{-1}$ throughout the $k$-th subproblem.

\begin{figure}
\centerline{\includegraphics[height=4cm,width=12cm]{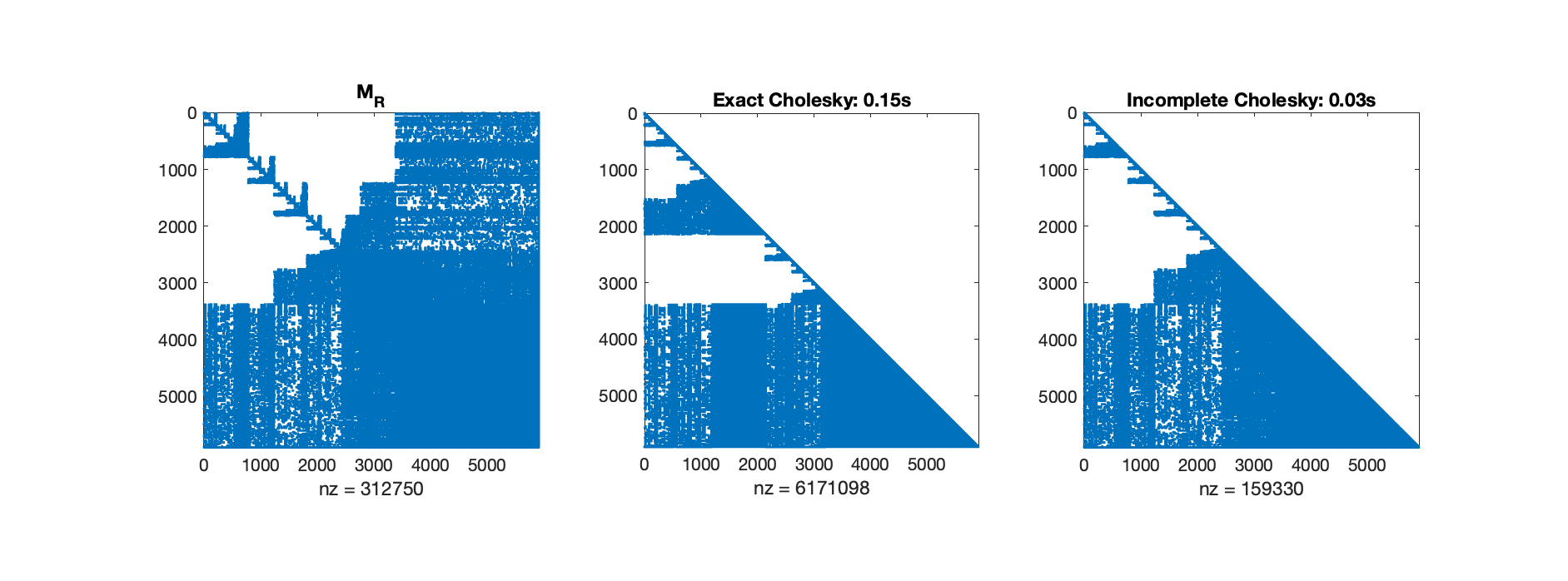}}
\caption{Cholesky factor of $M_R$ for the Lov\'asz Theta problem G51. PALM with exact Cholesky decomposition uses 28.4s, PALM with incomplete Cholesky decomposition uses 7.49s. } 
\label{fig:chol}
\end{figure}

When the weight matrix is chosen via the Cholesky decomposition of $M_R$, we may face high computational complexity and memory usage, as discussed in the previous subsection.  
To address this, we may instead use an incomplete Cholesky decomposition, for which our numerical results indicate is often an effective preconditioner in practice.  
In fact, the incomplete Cholesky decomposition is sometimes even better than the exact version, as it produces a sparse Cholesky factor $L$, thereby reducing the cost of forward and backward substitutions in computing $L^{-\top} L^{-1} x$. For example, for the Lov\'asz theta SDP instance on graph $G51$ reported in 
Table~\ref{test-illc}, the sparsity pattern of $M_R$ together with its exact 
and incomplete Cholesky factors is illustrated in Figure~\ref{fig:chol}. 
The incomplete Cholesky decomposition produces a significantly sparser factor 
and requires less computational time to construct. 
In our numerical test for this instance, PALM equipped with incomplete 
Cholesky preconditioning is approximately three times faster than that with the 
exact Cholesky decomposition.

In our implementation, the weight matrix in PALM is updated at each outer iteration. To compute it, we first perform a symbolic Cholesky factorization of $M_R$ using the {\sc Matlab} built-in function \texttt{symbfact} to estimate the number of nonzeros in its Cholesky factor. If this number is below $10^8$, we compute an exact Cholesky factorization; otherwise, we apply an incomplete Cholesky factorization via \texttt{ichol}. We simply set the tolerance for \texttt{ichol} to be zero, preventing any fill-ins in the factorization process.  Prior to calling \texttt{ichol}, we use \texttt{amd} to obtain a minimum-degree permutation, which reduces fill-in and alleviates the inaccuracy caused by incompleteness. To enhance numerical robustness, we add a small diagonal perturbation to $M_R$, namely $10^{-6}$ times its maximum diagonal entry, to ensure positive definiteness. The computational overhead of (incomplete) Cholesky factorization is typically modest due to the sparsity of $M_R$. Once the factorization is formed, evaluating the function value and gradient in the PALM subproblem only requires forward--backward substitutions, whose cost is significantly lower than that of computing the Cholesky factor itself.

It is worth noting that for certain SDP problems with an extremely large number of constraints $m$, forming the full matrix $M_R$ is already computationally prohibitive. Our current framework does not address this setting, and we therefore do not apply preconditioning in such cases. Instead, we avoid forming $M_R$ whenever its estimated number of nonzeros
exceeds $10^9$. Designing specialized preconditioners that 
exploit the structures in the constraint matrices and
the Gram matrix structure of $M_R$ without
the explicit construction of $M_R$ is a promising direction for future work. 

Another point is that, in the current implementation, we update the preconditioner at every outer iteration. However, if two consecutive iterates $R^k$ and $R^{k+1}$ are close, then the corresponding weight matrices $M_{R^k}$ and $M_{R^{k+1}}$ will also be close. In this case, one may reuse the previous preconditioner instead of recomputing the Cholesky factorization. This is another promising direction for future investigation.

\begin{remark}
Although the preconditioners in SDPF and PALM share the same matrix form, namely, the (exact or approximate) inverse of $M_{R^k},$ they are used in fundamentally different ways. In SDPF, the preconditioner is embedded in the PCG method used to solve the linear systems arising from the projection and retraction steps. In contrast, in PALM it is incorporated as a weight matrix in the penalty term. An important numerical observation is that SDPF typically requires an exact Cholesky factorization to keep the PCG iteration count below a prescribed threshold. In PALM, however, an incomplete Cholesky factorization is usually sufficiently effective and often lead to shorter run time, as discussed earlier. This makes PALM considerably more flexible in constructing the preconditioner.
\end{remark}

\subsection{Extension to general SDP problem}\label{Sec:gene}
In this subsection, we consider following more general SDP problem than 
\eqref{SDP}:
\begin{equation}\label{SDPG}
\min\left\{ f(X):\ \A(X)=b,\ \B(X)\in \P,\ X\in \S^n_+\right\}, 
\tag{SDPG}
\end{equation}
where $\B:\S^n\rightarrow \R^p$ is a linear map, $\P\subset \R^p$ is a closed convex set. By appropriately choosing the convex set $\P$, we can incorporate various constraints, such as inequality constraints and second-order cone constraints. The problem (\ref{SDPG}) can be reformulated as follows:
\begin{equation}\label{SDPG1}
\min\left\{ f(X)+\delta_{\P}(y):\ \A(X)=b,\ \B(X)=y,\ X\in \S^n_+,\ y\in \R^p\right\}, 
\end{equation}
where $\delta_{\P}(\cdot)$ is the indicator function of $\P,$ such that $\delta_{\P}(y)=0$ if $y\in \P$ and $\delta_{\P}(\cdot)=+\infty$ if $y\notin \P.$ Let $\mu\in \R^p$ be the Lagrange multipliers of the constraint $\B(X)=y$ in (\ref{SDPG}). The ALM subproblem of (\ref{SDPG1}) is given as follows:
\begin{equation}\label{ALMsubYZ}
\min\left\{f(X)+\delta_{\P}(y)-\<\mu,\B(X)-y\>+\frac{\beta}{2}\|\B(X)-y\|^2:\A(X)=b,X\in \S^n_+,\ y\in \R^p\right\}.
\end{equation}
where we preserve the constraints $\A(X)=b,X\in \S^n_+$ in the subproblem and penalize the other constraint $\B(X)-y$ with penalty parameters $\beta>0.$ When $X$ is fixed, the optimization problem (\ref{ALMsubYZ}) with respect to $y$ has the following closed form solution:
\begin{equation}\label{ALMsubYZcl}
y:=\Pi_\P\( \B(X)-\mu/\beta \):=\arg\min\left\{ \frac{1}{2}\|\B(X)-y-\mu/\beta\|^2:\ y\in \P \right\},
\end{equation}
where $\Pi_{\mathcal{P}}(\cdot)$ is the projection operator onto $\P.$ Substituting (\ref{ALMsubYZcl}) into (\ref{ALMsubYZ}), we get the following simplified ALM subproblem.
\begin{equation}\label{ALMsub}
\min\left\{f(X)+\frac{\beta}{2}\dist\(\B(X)-\mu/\beta,\P\)^2-\frac{\|\mu\|^2}{2\beta}:\A(X)=b,X\in \S^n_+\right\}.
\end{equation}
The objective function in (\ref{ALMsub}) is continuously differentiable with the following gradient \cite[Prop 12.30]{bauschke2017correction} 
\begin{equation}\label{ALMsubgrad}
\nabla f(X)+\beta\B^*\( \B(X)-\mu/\beta-\Pi_{\P}\(\B(X)-\mu/\beta\) \).
\end{equation}
Therefore, we can directly use a feasible method to solve (\ref{ALMsub}) to get the optimal solution $X^+$ and update the Lagrange multiplier as follows:
\begin{equation}\label{updatemu}
\mu^+ := \beta\(\Pi_\P\(\B(X)-\mu/\beta\)-\(\B(X)-\mu/\beta\)\).
\end{equation}
Apart from SDPF, Algorithm~\ref{alg:PALM} can also be directly extended to solve~\eqref{SDPG}.  
In Algorithm~\ref{alg:PALMG}, the preconditioner is applied only to the equality constraints, while the traditional penalty function is used for the remaining constraints. This design choice is motivated by the fact that our preconditioner is derived from the feasible method, which addresses only equality constraints. Developing suitable preconditioning strategies for inequality constraints remains an interesting topic for future research.

\begin{algorithm}
	\renewcommand{\algorithmicrequire}{\textbf{Input:}}
	\renewcommand{\algorithmicensure}{\textbf{Output:}}
	\caption{Preconditioned ALM for (\ref{SDPG})}
	\label{alg:PALMG}
	\begin{algorithmic}
		\STATE {\bf Initialization}: Choose $R^1\in \R^{n\times r},\lambda^1\in \R^m,\ \mu^1\in \R^p,\ \beta_1,\beta_2.$ 
		\FOR{$k = 1, 2, \dots$}
		\STATE 1. Choose a positive definite weight matrix $W^k\in \S^m.$ ($W^k\approx M_{R^k}^{-1}$ ideally).
		\STATE 2.	Solve the following ALM subproblem to get $R^{k+1}$
		\begin{align}
		\min\Big\{& f(RR^\top)-
        \<\lambda^k,\A(RR^\top)-b\>+\frac{\beta_1}{2}\|\A(RR^\top)-b\|^2_{W^k} \notag\\
		&+\frac{\beta_2}{2}\dist\(\B(RR^\top)-\mu^k/\beta_2,\P\)^2-\frac{\|\mu^k\|^2}{2\beta_2}:\ R\in \R^{n\times r} \Big\}. \label{SDPGALMsub}
		\end{align}
		\STATE 3. \quad $\lambda^{k+1}:=\lambda^k-\beta_1 W^k\(\A(R^{k+1}(R^{k+1})^\top)-b\).$	
		\begin{equation}
		\mu^{k+1}:=\beta_2\(\Pi_\P\(\B(R^{k+1}(R^{k+1})^\top)-\mu^k/\beta_2\)-\(\B(R^{k+1}(R^{k+1})^\top)-\mu^k/\beta_2\)\).  \notag 
		\end{equation}		
		\ENDFOR
	\end{algorithmic}
\end{algorithm}

\subsection{Convergence analysis}

In this subsection, we analyze the convergence of Algorithm~\ref{alg:PALMG}.  
The convergence of ALM for convex programming has been studied extensively~\cite{rockafellar1976augmented,cui2019r,xu2021iteration,wang2024local}, including cases where the subproblems are solved inexactly under various stopping criteria.  
However, to the best of our knowledge, there are no explicit results on ALM directly addressing a dynamically weighted quadratic penalty function.

Let $z=(X,y)\in \S^n\times \R^p=:\X,$ $\gamma:=(\lambda,\mu)\in \R^m\times \R^p$ and $\Q:=\S^n_+\times \P.$ We use $z^k:=(X^k,y^k),\gamma^k:=(\lambda^k,\mu^k)$ to denote the primal and dual variables in the $k$th iteration of Algorithm~\ref{alg:PALMG}. Define the following convex proper function
\begin{equation}\label{defig}
g(z):=f(X)+\delta_{\S^n_+}(X)+\delta_\P(y)=f(X)+\delta_\Q(z). 
\end{equation}
Let $\H:\S^n\times \R^p\rightarrow \R^m\times \R^p=:\Y$ be a linear operator such that
\begin{equation}\label{defiH} 
\H(z):=\( \A(X), \B(X)-y \). 
\end{equation}
Let $c:=(b,{\rm 0}_p)\in \Y.$ Problem (\ref{SDPG1}) can be rewritten as follows:
\begin{equation}\label{SDPG2}
\min\left\{ g(z):\ \H(z)=c,\ z\in \X  \right\}.
\end{equation}
Our convergence result is based on the following assumptions:
\begin{assumption}\label{assFeas}
The problem (\ref{SDPG2}) has a feasible solution $z^*\in \Q$ satisfying 
$\H(z^*)=c.$
\end{assumption}

\begin{assumption}\label{assALMsub}
For any $k\in \N^+,$ let $X^k:=R^kR^{k\top}.$ The $k$-th subproblem of Algorithm~\ref{alg:PALMG} is solved to sufficient accuracy such that  
\begin{equation}\label{optALMsub}
\dist\(\substack{\partial_X \delta_{\S^n_+}(X^{k+1})+\nabla f(X^{k+1})-\A^*(\lambda^k)+\beta_1 W_k\( \A(X^{k+1})-b \)\\
+\beta_2\B^*\( \B(X^{k+1})-\mu^k/\beta_2-\Pi_{\P}\(\B(X^{k+1})-\mu^k/\beta_2\) \)},0\)<\epsilon_k
\end{equation}
for some $\epsilon_k>0$ such that $\sum_{k=1}^\infty\epsilon_k<\infty.$ 
\end{assumption}
The condition~\eqref{optALMsub} in Assumption~\ref{assALMsub} is the approximate stationarity condition of the ALM subproblem~\eqref{SDPGALMsub}, with the term $RR^\top$ replaced by $X \in \mathbb{S}^n_+$. Such a point can be obtained using a rank-adaptive first-order method~\cite{tang2024feasible}, where the rank parameter $r$ is dynamically updated to escape saddle points. For simplicity, we assume in Algorithm~\ref{alg:PALMG} that $r \in \mathbb{N}^+$ is fixed. 

Now we state and prove the convergence theorem of Algorithm~\ref{alg:PALMG}. Our proof is a direct adaptation of the convergence analysis of \cite{rockafellar1976augmented} by Rockafellar. For completeness, we include the proof in Section~\ref{Proof-of-theorem-1}.

\begin{theorem}\label{convtheo}
Suppose Assumption~\ref{assFeas} and Assumption~\ref{assALMsub} holds. If the sequence $(z^k,\gamma^k)$ generated by Algorithm~\ref{alg:PALMG} satisfies that $\{z^k\}_{k\in \N^+}$ is bounded and there exists $\sigma > 0$ such that $W^k\succeq \sigma I_m$ for all $k$, then the following KKT residual satisfies 
\begin{equation}\label{ALMKKTres}
\max\left\{ \dist\( \partial_{z} g(z^k)-\H^*(\gamma^k),0 \),\ \|\H(z^k)-c\|_2 \right\} \; \to \; 0, \quad \mbox{as $k\to \infty$.}
\end{equation}
\end{theorem}
\section{Solver description of SDPF+}\label{Sec-solver}

The previous section reviewed the feasible method SDPF and introduced the PALM. In this section, we present the new {\sc Matlab}  solver SDPF+, which employs these two algorithms in consecutive stages.

For any given SDP problem, the algorithm initially applies the first stage -- the feasible method SDPF -- until one of the following conditions is encountered:
\begin{itemize}
\item[(1)] When exact Cholesky decomposition is too costly, the preconditioned conjugate gradient method for solving~\eqref{linsolve} may not reach the required accuracy within the iteration limit, making the projection and retraction steps unreliable. In our implementation, we perform a symbolic Cholesky factorization to estimate the number of nonzeros in the Cholesky factor. If this number exceeds the prescribed threshold $10^8,$ we deem the factorization too expensive. In this situation, only a diagonal preconditioner is used in the PCG during the SDPF stage. If PCG does not converge within 20 iterations, we declare it unsuccessful and switch to PALM with an incomplete Cholesky preconditioner.

\item[(2)] In~\cite{tang2024feasible}, the retraction is computed iteratively using the Gauss--Newton method. Although it exhibits local quadratic convergence, convergence within the preset maximum number of iterations, which we set to be 20, is not guaranteed. If the Gauss--Newton method does not converge after 20 steps, we switch to PALM.

\item[(3)] For degenerate SDP problems, the feasible method may fail to achieve the required dual feasibility~\cite[Section~5]{tang2024feasible}. This is caused by the inaccuracy in solving the nearly singular linear system~\eqref{linsolve}, which in turn degrades the accuracy of the recovered Lagrange multipliers. In our implementation, we declare an SDP instance as (approximately) degenerate if during the SDPF stage, PCG with a diagonal preconditioner fails to terminate within 20 iterations.
\end{itemize}

If any of the above conditions occur, it is futile to use the feasible method SDPF, and the algorithm switches to the second stage, which employs the preconditioned augmented Lagrangian method (PALM). The complete algorithmic framework of SDPF+  is presented in Algorithm~\ref{alg:SDPF+}.

\begin{algorithm}
	\renewcommand{\algorithmicrequire}{\textbf{Input:}}
	\renewcommand{\algorithmicensure}{\textbf{Output:}}
	\caption{SDPF+}
	\label{alg:SDPF+}
	\begin{algorithmic}
		\STATE {\bf Input}:  problem (\ref{SDPG}). 
		\IF{There is no constraint $X\in \P$ or $\H(X)\in \Q$}
		\STATE Use the feasible method SDPF to solve the problem (\ref{SDP}) and switch to Algorithm~\ref{alg:PALM} when any of the three conditions mentioned above occurs.
		\ELSE 
		\STATE Use ALM to solve the problem (\ref{SDPG}), with the subproblem (\ref{ALMsub}) solved by the feasible method SDPF and switch to Algorithm~\ref{alg:PALMG} when any of the above three conditions mentioned above occurs.
		\ENDIF
	\end{algorithmic}
\end{algorithm}

In Algorithm~\eqref{alg:SDPF+}, we retain SDPF as the first stage solver rather than directly applying PALM because, for certain well-conditioned SDP problems with block-diagonal constraints, such as Max-Cut and rotation synchronization, PCG with a diagonal preconditioner (or even without preconditioning) converges within fewer than 10 iterations. For these cases, SDPF is particularly efficient, as it avoids forming $M_R$ and computing its Cholesky factorization. Actually, for well-conditioned block-diagonal SDP problems, both PALM and SDPF are already efficient even without preconditioning. However, SDPF still enjoys several advantages over PALM due to its single-loop and penalty-free nature, which makes the algorithm more robust. In contrast, the performance of PALM can be sensitive to the update of the penalty parameter and the stopping tolerances of the subproblems.

Let the dual slack matrix be 
\begin{equation}\label{dSlack}
S:=\nabla f(X)-\A^* (\lambda)-\B^*(\mu).
\end{equation}
We use the following KKT residuals to measure the accuracy of the 
iterate $(X,\lambda,\mu)$ for (\ref{SDPG}):
 \begin{align}
 {\rm pfeas} &:=\max\left\{\frac{\|\A(X)-b\|_2}{1+\|b\|_2},\ \frac{\|\B(X)-y\|_2}{1+\|\B(X)\|_2}\right\},
 \label{pfeas} 
 \\
 {\rm dfeas} &:=\frac{\|\Pi_{\S^n_+}(-S)\|_F}{1+\|S\|_F},\quad {\rm comp}:=\frac{|\<X,S\>|}{1+\|\nabla f(X)\|_F}, 
 \label{dfeascomp}
 \end{align}
 where $y$ is given in (\ref{ALMsubYZcl}).
 We terminate the algorithm when the KKT residuals fall below a specified tolerance. Note that we omit checking the residual for the constraint $X \in \mathbb{S}^n_+$, $y\in \P$ as well as the dual residual and complementarity conditions associated with the variable $y$, since these conditions are satisfied exactly at every iteration of SDPF+.

When the problem \eqref{SDPG} is a standard linear SDP problem, that is, $f(X) = \langle C, X \rangle$ for some $C \in \mathbb{S}^n$ and there are no constraints $\B(X)\in \P$, we also compute the following duality gap:
\begin{equation}\label{pdgap}
{\rm pdgap}:=\frac{|\lambda^\top b-\<C,X\>|}{1+|\lambda^\top b|+|\<C,X\>|}
\end{equation} 
 and use it to replace the complementarity gap when ${\rm pdgap}<{\rm comp}.$
 
We refer the reader to Section~\ref{App-solver}
for the usage and data structures of our {\sc Matlab} implementation of SDPF+,
where it can also solve the multi-block version of \eqref{SDPG}.

\section{Numerical experiments}\label{Sec-numer}

In this section, we conduct various numerical experiments to verify the efficiency and robustness of SDPF+. For linear SDP problems, the solvers we compare with are:
\begin{itemize}
\item SDPNAL+ \cite{SDPNALp2}: a linear SDP solver based on dual ADMM and dual ALM. 
\item SDPT3 \cite{TTT}: A linear SDP solver based on primal-dual interior point methods. 
\item ManiSDP \cite{wang2025solving}: A linear SDP solver based on low-rank ALM. 
\end{itemize}
While both SDPF+ and ManiSDP utilize low-rank ALM, the main difference lies in 
the use of the  penalty function: SDPF+ employs a preconditioner weighted penalty function, whereas ManiSDP adopts the traditional penalty function. By comparing SDPF+ with ManiSDP, we aim to assess the effectiveness of preconditioning strategy in low-rank ALM.

For nonlinear SDP problems, we compare SDPF+ with specialized solvers or algorithms tailored to the specific problem instances. These algorithms will be described in the subsequent subsections. Unless otherwise specified, the maximum running time is set to 3600 seconds, and the termination tolerance for the maximum KKT residual is $10^{-6}$. All experiments are conducted using {\sc Matlab} R2021b on a workstation equipped with an Intel Xeon E5-2680 v3 @ 2.50GHz processor and 128GB RAM.

For each solver and problem instance, we report the KKT residuals, function values, and run time (wall clock times). If a solver fails to achieve the maximum KKT residual of less than $10^{-2}$, we omit its output. We use bold fonts to indicate the most efficient algorithm achieving the required accuracy for each instance. For problem classes with many test instances, we additionally use performance profiles to visually compare the efficiency and accuracy of different algorithms.

\subsection{Lov\'asz Theta problem}

We test the following Lov\'asz Theta problem, which is frequently used to benchmark SDP solvers:
\begin{equation}\label{ThetaSDP}
\min\left\{ -\<{\bf 1}_{n\times n},X\>:\ \<I,X\>=1,\ X_{ij}=0,\ \forall ij\in \E_G,\ X\in \S^n_+ \right\},
\end{equation}
where $\E_G$ is the edge set of a graph $G$ with $n$ vertices. The dataset we test comes from coding theory \cite{Coding}. 
 
\begin{center}
\begin{tiny}
\begin{longtable}{|c|ccccc|c|}
\caption{Test results for Lov\'asz Theta SDP problems.}
\label{test-theta}
\\
\hline
problem & algorithm & pfeas& dfeas & comp & fval & time \\ \hline
\endhead
1dc.1024  & {\bf SDPF+} & 6.78e-07 & 4.74e-07 & 9.90e-08 & -9.5986179e+01& 4.62e+01 \\
n=1024& SDPNAL+ & 9.52e-07 & 7.33e-08 & 8.52e-09 & -9.5986460e+01& 9.67e+01 \\
m=24064& SDPT3 & 5.86e-12 & 0.00e+00 & 8.67e-08 & -9.5984581e+01& 6.25e+02 \\
rank=163 & ManiSDP & 9.95e-07 & 2.35e-07 & 3.30e-18 & -9.5987258e+01& 5.84e+02 \\
 \hline
1dc.2048 & {\bf SDPF+} & 5.33e-07 & 8.42e-07 & 2.07e-07 & -1.7473140e+02& 1.94e+02 \\
n=2048 & SDPNAL+ & 3.12e-07 & 7.72e-08 & 1.03e-09 & -1.7473104e+02& 9.32e+02 \\
m=58368&SDPT3 & - & - & - & -& - \\
rank=256 & ManiSDP & 3.36e-06 & 2.93e-07 & 1.68e-17 & -1.7475532e+02& 3.61e+03 \\
\hline
1et.1024  & SDPF+ & 6.65e-07 & 7.99e-07 & 9.48e-14 & -1.8422910e+02& 1.04e+02 \\
n=1024    & SDPNAL+  & 6.28e-07 & 6.36e-09 & 2.64e-09 & -1.8422783e+02 & 3.71e+02 \\
m=9601    & {\bf SDPT3}    & 1.47e-12 & 0.00e+00 & 1.06e-07 & -1.8422591e+02 & 7.86e+01 \\
rank=226 & ManiSDP & 1.68e-06 & 4.66e-05 & 4.85e-17 & -1.8422889e+02& 3.61e+03 \\
\hline
1et.2048  & {\bf SDPF+} & 9.98e-07 & 3.07e-07 & 3.20e-11 & -3.4203828e+02& 2.16e+02 \\
n=2048    & SDPNAL+  & 6.51e-07 & 1.18e-08 & 1.04e-08 & -3.4203305e+02 & 1.91e+03 \\
m=22529   & SDPT3    & 4.74e-13 & 0.00e+00 & 2.37e-07 & -3.4202835e+02 & 5.50e+02 \\
rank=296& ManiSDP & 2.21e-05 & 4.23e-06 & 2.97e-17 & -3.4210737e+02& 3.60e+03 \\
\hline
1tc.1024  & SDPF+ & 9.89e-07 & 7.92e-07 & 3.35e-10 & -2.0630739e+02& 1.29e+02 \\
n=1024    & SDPNAL+  & 6.34e-07 & 8.30e-09 & 4.54e-06 & -2.0631074e+02 & 5.37e+02 \\
m=7937    & {\bf SDPT3}    & 2.75e-14 & 0.00e+00 & 2.08e-07 & -2.0630400e+02 & 5.30e+01 \\
rank=200 & ManiSDP & 4.04e-06 & 7.34e-07 & 1.02e-16 & -2.0632068e+02& 3.60e+03 \\
\hline
1tc.2048 & {\bf SDPF+} & 9.94e-07 & 2.37e-07 & 1.12e-14 & -3.7465327e+02& 2.54e+02 \\
n=2048    & SDPNAL+  & 4.69e-05 & 2.71e-07 & 1.44e-09 & -3.7465245e+02 & 3.64e+03 \\
m=18945   & SDPT3    & 9.22e-13 & 0.00e+00 & 1.61e-07 & -3.7464286e+02 & 3.96e+02 \\
rank=281 & ManiSDP & 1.28e-05 & 2.15e-06 & 2.44e-17 & -3.7482453e+02& 3.61e+03 \\
\hline
1zc.1024  & SDPF+ & 1.07e-07 & 5.81e-07 & 4.11e-12 & -1.2866665e+02& 5.40e+01 \\
n=1024    & {\bf SDPNAL+}  & 2.37e-07 & 5.06e-08 & 1.48e-08 & -1.2866636e+02 & 4.23e+01 \\
m=16641   & SDPT3    & 2.28e-15 & 0.00e+00 & 1.24e-07 & -1.2866657e+02 & 1.83e+02 \\
rank=294 & ManiSDP & 6.08e-07 & 1.85e-07 & 2.02e-17 & -1.2866718e+02& 2.69e+02 \\
\hline
1zc.2048  & {\bf SDPF+} & 1.67e-07 & 4.33e-07 & 8.93e-08 & -2.3740040e+02& 6.29e+01 \\
n=2048    & SDPNAL+  & 6.05e-07 & 1.59e-09 & 6.11e-11 & -2.3740044e+02 & 4.06e+02 \\
m=39425   & SDPT3    & -        & -        & -        & -              & -        \\
rank=396& ManiSDP & 8.65e-07 & 3.08e-07 & 1.58e-18 & -2.3740255e+02& 4.72e+02 \\
\hline
1zc.4096  & {\bf SDPF+} & 4.67e-07 & 9.85e-07 & 2.67e-09 & -4.4916594e+02& 1.19e+03 \\
n=4096    & SDPNAL+  & 6.55e-07 & 1.73e-08 & 3.70e-09 & -4.4916755e+02 & 2.91e+03 \\
m=92161   & SDPT3    & -        & -        & -        & -              & -        \\
rank=799 & ManiSDP & 2.40e-05 & 1.11e-06 & 2.52e-17 & -4.4944055e+02& 3.64e+03 \\
\hline
\end{longtable}
\end{tiny}
\end{center}

From Table~\ref{test-theta}, we observe that SDPF+ successfully solves all instances to the required accuracy and is the most efficient solver for 6 out of the 9 problems. In many cases, SDPF+ is more than 10 times faster than ManiSDP, demonstrating the effectiveness of the preconditioner in low-rank ALM.

We also test the following Theta+ problem
\begin{equation}\label{ThetaDNN}
\min\left\{ -\<{\bf 1}_{n\times n},X\>:\ \<I,X\>=1,\ X_{ij}=0,\ \forall ij\in \E_G,\ X\in \S^n_+,\ X\geq 0 \right\},
\end{equation}
which is obtained by adding the nonnegativity constraint $X\in\P = \{ X\in\S^n : X\geq 0\}$ to (\ref{ThetaSDP}). Because ManiSDP cannot handle inequality constraints and SDPT3 based on interior point method is extremely 
expensive for problems with $\Omega(n^2)$ constraints, we only compare SDPF+ with SDPNAL+.

\begin{center}
\begin{tiny}
\begin{longtable}{|c|ccccc|c|}
\caption{Test results for Theta+ problem.}
\label{test-theta+}
\\
\hline
problem & algorithm & pfeas& dfeas & comp & fval & time \\ \hline
\endhead
1dc.1024  & SDPF+ & 7.89e-07 & 9.17e-07 & 4.28e-11 & -9.5552604e+01& 4.46e+02 \\
rank=386    & {\bf SDPNAL+}  & 5.42e-07 & 2.84e-07 & 2.92e-08 & -9.5551160e+01 & 1.76e+02 \\
\hline
1dc.2048  & SDPF+    & 7.01e-07 & 7.51e-07 & 4.99e-10 & -1.7426269e+02& 2.38e+03 \\
rank=670    & {\bf SDPNAL+}  & 9.70e-07 & 5.37e-07 & 8.77e-07 & -1.7425908e+02 & 1.32e+03 \\
\hline
1et.1024  & SDPF+     & 5.67e-07 & 8.72e-07 & 1.08e-08 & -1.8207581e+02& 3.96e+02 \\
rank=365    & {\bf SDPNAL+}  & 9.33e-07 & 4.62e-07 & 3.66e-06 & -1.8207513e+02 & 1.44e+02 \\
\hline
1et.2048  & SDPF+    & 1.77e-06 & 6.05e-07 & 2.30e-08 & -3.3819558e+02& 3.60e+03 \\
rank=640    & {\bf SDPNAL+}  & 1.31e-06 & 4.68e-07 & 4.57e-06 & -3.3817316e+02 & 1.64e+03 \\
\hline
1tc.1024  & SDPF+     & 1.00e-06 & 4.63e-07 & 2.43e-09 & -2.0420941e+02& 1.10e+03 \\
rank=346    & {\bf SDPNAL+}  & 1.40e-06 & 5.90e-07 & 5.51e-06 & -2.0420916e+02 & 4.00e+02 \\
\hline
1tc.2048  & SDPF+     & 4.50e-07 & 9.78e-07 & 3.34e-10 & -3.7049547e+02& 2.69e+03 \\
rank=603    & {\bf SDPNAL+}  & 9.11e-07 & 6.90e-07 & 1.22e-06 & -3.7049077e+02 & 2.05e+03 \\
\hline
1zc.1024  & SDPF+     & 9.38e-07 & 3.97e-07 & 3.07e-09 & -1.2800306e+02& 3.51e+02 \\
rank=584    & {\bf SDPNAL+}  & 9.56e-07 & 8.08e-07 & 2.56e-06 & -1.2799771e+02 & 5.31e+01 \\
\hline
1zc.2048  & SDPF+    & 1.98e-07 & 9.94e-07 & 5.04e-10 & -2.3740062e+02& 5.23e+02 \\
rank=505    & {\bf SDPNAL+}  & 5.07e-07 & 2.36e-07 & 5.68e-11 & -2.3739990e+02 & 1.45e+02 \\
\hline
1zc.4096  & SDPF+     & 7.07e-06 & 4.34e-06 & 2.60e-07 & -4.4377186e+02& 3.60e+03 \\
rank=1326    & SDPNAL+  & 7.80e-06 & 5.22e-06 & 4.74e-08 & -4.4373543e+02 & 3.61e+03 \\
\hline
\end{longtable}
\end{tiny}
\end{center}
From Table~\ref{test-theta+}, we observe that SDPF+ is slower than SDPNAL+. This result is expected, as the presence of nonnegativity constraints usually leads to higher solution ranks compared to the standard SDP formulation~(\ref{ThetaSDP}). In some instances, the rank reaches nearly one-third, and occasionally even one-half, of the matrix dimension, thus posing a challenge for low-rank SDP solvers. Moreover, SDPNAL+ is well known for its efficiency in solving doubly nonnegative SDP problems. Nonetheless, the performance of SDPF+ remains satisfactory, as its run time is typically within three times that of SDPNAL+, demonstrating its robustness across challenging problem instances.

\subsection{Quadratic assignment problem}

We test the doubly nonnegative relaxation of the quadratic assignment problem (QAP):
\begin{equation}\label{QAP}
\min\left\{ \<X,AXB\>:\ X^\top X=I_n,\ X\in \{0,1\}^{n\times n} \right\}.
\end{equation}
There exist several formulations of doubly nonnegative (DNN) relaxations for the QAP, which have been shown to be equivalent~\cite{povh2009copositive}. For SDPNAL+, we use the $QAP_{R_3}$ relaxation described in~\cite[Section~4]{povh2009copositive}, as this formulation has been adopted in previous numerical experiments by SDPNAL+~\cite{SDPNALp2, SDPNALp1, SDPNAL,QSDPNAL} and is known for its superior performance. For SDPF+, we choose the equivalent $QAP_{LS}$ relaxation from~\cite[Section~5]{povh2009copositive}, which includes a facial constraint of the form $XU = 0$ \cite{graham2022restricted}. As discussed in Appendix~\ref{App-solver}, SDPF+ can exploit this structure to handle the constraint efficiently under the low-rank factorization setting.
We test all instances in the QAPLIB benchmark~\cite{burkard1997qaplib} with problem sizes satisfying $32 \leq n \leq 40$.

\begin{center}
\begin{tiny}
\begin{longtable}{|c|ccccc|c|}
\caption{Test results for convex relaxation of quadratic assignment problem (\ref{QAP}). The coefficient matrices $A,B$ are normalized by their maximum entries.}
\label{test-QAP}
\\
\hline
problem & algorithm & pfeas& dfeas & comp & fval & time \\ \hline
\endhead
 esc32a    & SDPF+ & 3.31e-07 & 9.71e-07 & 1.03e-09 & 6.7789147e-02& 2.52e+03 \\
  rank=827  & {\bf SDPNAL+}  & 1.01e-06 & 1.92e-07 & 1.47e-06 & 6.7789339e-02 & 2.49e+02 \\
\hline
esc32b    & SDPF+ & 6.29e-07 & 9.18e-07 & 1.44e-08 & 1.0646688e-01& 5.51e+02 \\
  rank=730  & {\bf SDPNAL+}  & 1.00e-06 & 1.74e-07 & 2.24e-06 & 1.0642686e-01 & 2.74e+02 \\
\hline
esc32c   & {\bf SDPF+} & 4.41e-07 & 6.97e-07 & 1.03e-08 & 2.5019981e-01& 1.20e+02 \\
  rank=565 & SDPNAL+  & 7.40e-07 & 4.27e-07 & 6.82e-06 & 2.5018622e-01 & 2.16e+02 \\
\hline
esc32d    & SDPF+ & 5.62e-07 & 7.91e-07 & 4.57e-09 & 1.7517284e-01& 3.55e+02 \\
  rank=660  & {\bf SDPNAL+}  & 6.14e-07 & 3.06e-07 & 3.02e-06 & 1.7517317e-01 & 4.64e+01 \\
\hline
esc32e    & SDPF+ & 8.70e-08 & 6.87e-07 & 5.46e-08 & 2.9102252e-03& 8.65e+01 \\
  rank=540 & {\bf SDPNAL+}  & 2.38e-07 & 5.60e-07 & 1.96e-08 & 2.9101854e-03 & 8.33e+01 \\
\hline
esc32g  & SDPF+ & 3.34e-07 & 6.32e-07 & 9.23e-10 & 1.3079328e-02& 1.08e+02 \\
rank=560 & {\bf SDPNAL+} & 4.24e-07 & 5.88e-07 & 1.82e-08 & 1.3079318e-02& 4.12e+01 \\
\hline
esc32h & SDPF+ & 6.63e-07 & 9.37e-07 & 3.72e-08 & 2.5145117e-01& 1.01e+03 \\
rank=843 & {\bf SDPNAL+} & 8.34e-07 & 4.04e-07 & 2.25e-05 & 2.5139278e-01& 8.23e+02 \\
\hline
kra32       & SDPF+ & 4.48e-07 & 9.73e-07 & 2.39e-09 & 3.0164203e-01& 1.61e+03 \\
  rank=596 & {\bf SDPNAL+}  & 1.21e-06 & 2.01e-07 & 2.05e-05 & 3.0147487e-01 &  6.66e+02 \\
\hline
lipa40a    & SDPF+ & 6.42e-07 & 7.14e-07 & 5.82e-07 & 8.7819603e-01& 1.99e+03 \\
  rank=840 & {\bf SDPNAL+}  & 3.32e-07 & 1.65e-07 & 2.02e-09 & 8.7819411e-01 & 1.27e+03 \\
\hline
lipa40b   & SDPF+ & 6.99e-07 & 0.00e+00 & 1.53e-08 & 6.1038000e-01& 5.52e+02 \\
  rank=832 & SDPNAL+  & 1.17e-07 & 5.53e-07 & 3.60e-09 & 6.1037986e-01 & 5.28e+02 \\
\hline
sko42     & SDPF+ & 4.65e-06 & 6.86e-06 & 8.61e-07 & 4.7098467e-01& 3.60e+03 \\
 rank=1113 & {\bf SDPNAL+}  & 1.52e-06 & 3.63e-07 & 2.59e-06 & 4.7094929e-01 & 1.54e+03 \\
\hline
ste36a   & SDPF+ & 2.65e-05 & 1.18e-05 & 1.34e-07 & 8.8679914e-02& 3.60e+03 \\
 rank=826  & {\bf SDPNAL+}  & 1.60e-06 & 4.30e-07 & 3.09e-06 & 8.8669185e-02 & 1.06e+03 \\
\hline
ste36b    & SDPF+ & 1.81e-05 & 1.02e-05 & 7.35e-08 & 3.2160132e-02& 3.60e+03 \\
   rank=828  & {\bf SDPNAL+}  & 9.03e-07 & 5.12e-07 & 3.19e-06 & 3.2181769e-02 & 2.41e+03 \\
\hline
ste36c    & SDPF+ & 1.40e-05 & 7.08e-06 & 1.98e-07 & 9.4396324e-02& 3.60e+03 \\
  rank=828 & {\bf SDPNAL+}  & 1.56e-06 & 4.07e-07 & 3.30e-06 & 9.4389140e-02 & 1.33e+03 \\
\hline
tai35a    & SDPF+ & 6.69e-07 & 1.99e-07 & 6.88e-10 & 5.6398422e-01& 4.31e+02 \\
 rank=650  & {\bf SDPNAL+}  & 1.06e-06 & 1.62e-07 & 2.70e-05 & 5.6395534e-01 & 3.80e+02 \\
\hline
tai35b    & SDPF+ & 1.32e-05 & 1.98e-05 & 4.46e-07 & 1.3004610e-01& 3.60e+03 \\
  rank=833 & {\bf SDPNAL+}  & 6.82e-07 & 6.84e-07 & 7.30e-06 & 1.3000286e-01 & 2.71e+03 \\
\hline
tai40a     & SDPF+ & 9.71e-07 & 2.63e-07 & 8.26e-11 & 5.5824653e-01& 9.06e+02 \\
 rank=849 & {\bf SDPNAL+}  & 1.20e-06 & 1.48e-07 & 1.60e-05 & 5.5822841e-01 & 8.90e+02 \\
\hline
tai40b      & SDPF+ & 1.31e-05 & 5.33e-05 & 5.80e-07 & 1.3799361e-01& 3.60e+03 \\
  rank=1222 & {\bf SDPNAL+}  & 1.56e-06 & 6.94e-07 & 1.95e-05 & 1.3794748e-01 &  3.34e+03 \\
\hline
tho40     & SDPF+     & 3.49e-06 & 1.43e-06 & 9.34e-11 & 3.2487741e-01 & 3.60e+03 \\
  rank=946  & {\bf SDPNAL+}  & 8.91e-07 & 2.37e-07 & 1.76e-05 & 3.2480955e-01 & 1.32e+03 \\
\hline
\end{longtable}
\end{tiny}
\end{center}

From Table~\ref{test-QAP}, we observe that SDPF+ is slower than SDPNAL+, likely due to the fact that all the tested problems have high-rank solutions, with ranks exceeding half of the matrix dimension. To the best of our knowledge, such high-rank instances are rarely considered in the numerical experiments of low-rank SDP solvers. Nevertheless, the performance of SDPF+ is quite satisfactory: its run time is typically within three times that of SDPNAL+, and it often achieves higher solution accuracy. This is notable given that SDPNAL+ is widely regarded as a state-of-the-art solver for DNN relaxations of  QAP problems. 

In practice, for SDP problems with a large number of constraints—such as 
DNN relaxations of certain combinatorial optimization problems (e.g., QAP 
and the maximum stable set problem), the optimal solution may have high rank. 
In such cases, the low-rank solver SDPF+ can be less efficient than SDPNAL+, 
and we recommend using SDPNAL+. On the other hand, for SDP problems with relatively few constraints that 
guarantee the existence of low-rank optimal solutions \cite{BAI,pataki1998on}, SDPF+ is expected to 
perform well and is recommended. We note, however, that a large number of 
constraints does not necessarily destroy low-rank structure, as illustrated 
in the  two classes of examples
in subsections~\ref{subsec:graph-equipartition} and ~\ref{subsec:rank-1-tensor}, where SDPF+ remains efficient. In practice, predicting whether an SDP problem will admit a low-rank solution 
prior to solving it is generally difficult, as this requires a good 
understanding of the problem structure. Nevertheless, even for high-rank 
instances, SDPF+ remains reasonably robust in our experiments, with running 
times typically within a factor of five that of SDPNAL+.

\subsection{Graph equipartition}
\label{subsec:graph-equipartition}

In this subsection, we test the DNN relaxations of graph equipartition problems \cite{rendl1999semidefinite,sotirov2014efficient}:
\begin{equation}\label{GEP}
\min\left\{\<L_G,X\>:\ X{\bf 1}_n=\frac{n}{K}{\bf 1}_n,\ \dd(X)={\bf 1}_n,\ X\in \S^n_+,\ X\geq 0\right\},
\end{equation}
where $L_G$ is the Laplacian matrix of a graph and $K$ is the number of clusters. Because the constraint $X\mathbf{1}_n = \frac{n}{K}\mathbf{1}_n$ implies that $\mathbf{1}_n$ is an eigenvector of $X$, we can, following the procedure in~\cite{tang2024solving}, reformulate problem~(\ref{GEP}) into an equivalent but simpler problem with reduced number of constraints:
\begin{equation}\label{GEP1}
\min\left\{\<L_G,X\>:\ \<X,{\bf 1}_{n\times n}\>=0,\ \dd(X)={\bf 1}_n,\ X\in \S^n_+,\ X\geq \frac{-1}{K-1}\right\}.
\end{equation}
We randomly generated a graph with $K$ clusters of equal size, where the edge density within each cluster is $0.2$ and the edge density between different clusters is $0.1$. This graph generation procedure, commonly known as the stochastic block model, is frequently used in numerical simulations for graph partitioning~\cite{abbe2018community}.

\begin{center}
\begin{tiny}
\begin{longtable}{|c|ccccc|c|}
\caption{Test results for DNN relaxations of graph equipartition problems.}
\label{test-GEP}
\\
\hline
problem & algorithm & pfeas& dfeas & comp & fval & time \\ \hline
\endhead
n=300& {\bf SDPF+} & 3.97e-07 & 6.37e-07 & 3.20e-07 & 8.5204105e+03& 1.87e+01 \\
K=3,rank=35& SDPNAL+ & 2.51e-08 & 7.26e-07 & 7.25e-06 & 8.5204005e+03& 2.68e+01 \\
\hline
n=500& {\bf SDPF+} & 8.66e-07 & 5.69e-07 & 7.15e-08 & 2.3837915e+04& 3.61e+01 \\
 K=5,rank=89& SDPNAL+ & 3.65e-07 & 4.96e-07 & 3.25e-07 & 2.3837913e+04& 1.84e+02 \\
\hline
 n=1000 & {\bf SDPF+} & 7.38e-07 & 6.25e-07 & 3.70e-10 & 9.3643421e+04& 1.19e+02 \\
K=10,rank=228& SDPNAL+ & 3.40e-07 & 2.74e-07 & 1.93e-06 & 9.3643408e+04& 1.01e+03 \\
\hline
 n=600 & {\bf SDPF+} & 5.21e-07 & 5.12e-07 & 7.70e-09 & 3.5763287e+04& 1.16e+02 \\
K=3,rank=49& SDPNAL+ & 2.52e-09 & 1.46e-07 & 3.55e-06 & 3.5763274e+04& 1.22e+02 \\
\hline
 n=1000 & {\bf SDPF+} & 5.00e-07 & 3.64e-07 & 6.43e-07 & 9.8624195e+04& 2.45e+02 \\
K=5,rank=148& SDPNAL+ & 7.77e-07 & 5.26e-07 & 4.08e-06 & 9.8624177e+04& 3.27e+02 \\
\hline
 n=2000 & {\bf SDPF+} & 2.21e-07 & 6.90e-07 & 8.03e-09 & 3.9223713e+05& 5.65e+02 \\
K=10,rank=318& SDPNAL+ & 1.48e-07 & 4.05e-07 & 1.68e-06 & 3.9223710e+05& 3.33e+03 \\
\hline
 n=1500  & SDPF+ & 9.34e-07 & 0.00e+00 & 3.16e-08 & 2.2561800e+05& 3.59e+02 \\
K=3,rank=7& {\bf SDPNAL+} & 1.60e-07 & 2.76e-07 & 7.74e-07 & 2.2561797e+05& 2.81e+02 \\
\hline
n=2500 & {\bf SDPF+} & 3.67e-07 & 1.18e-09 & 9.22e-07 & 6.2291001e+05& 1.03e+03 \\
K=5,rank=7& SDPNAL+ & 2.27e-08 & 1.12e-05 & 1.15e-06 & 6.2291036e+05& 3.60e+03 \\
\hline
 n=5000& SDPF+ & 5.28e-05 & 8.41e-04 & 2.50e-04 & 2.4990046e+06& 3.61e+03 \\
K=10,rank=523& SDPNAL+ & 6.40e-05 & 3.61e-04 & 2.63e-04 & 2.5388198e+06& 3.60e+03 \\
\hline
\end{longtable}
\end{tiny}
\end{center}

From Table~\ref{test-GEP}, we observe that SDPF+ is usually more efficient than SDPNAL+. This improved efficiency is attributed to the fact that the solution ranks of problem~(\ref{GEP1}) are lower compared to the SDP problems with bound constraints that we have tested in the last subsection.

\subsection{Rank-1 tensor approximation}
\label{subsec:rank-1-tensor}

We test the SDP relaxation of the rank-1 tensor approximation problem \cite{nie2014semidefinite}: 
\begin{equation}\label{R1tensor}
\min\left\{ \<C,X\>:\ \A(X)=b,\ X\in \S^n_+ \right\}.
\end{equation}
In the SDP relaxation for the best rank-1 tensor approximation problem \cite{nie2014semidefinite}, 
the linear constraints $\mathcal{A}(X)=b$ arise from lifting the unit-sphere constraint 
and enforcing algebraic consistency among monomials in the moment matrix. 
The associated constraint matrices are highly structured and typically sparse, 
as each constraint involves only a few entries of $X$. 
Nevertheless, the total number of constraints is in the order of $\Omega(n^2)$, 
since nearly all entries of $X$ participate in consistency conditions.

We do not compare with SDPT3 because this problem has large number of constraints, which is too expensive to be solved by interior-point methods. Because the optimal solution of this problem usually has rank equal to $n-1,$ we use SDPF+ and ManiSDP to solve the following dual problem of (\ref{R1tensor}), which has a rank-1 optimal solution
\begin{equation}\label{R1tensord}
\max\left\{ b^\top y:\ S+\A^* y=C,\ S\in \S^n_+,\ y\in \R^m \right\}.
\end{equation}
To eliminate the free variable $y,$ let $\hat{\A}:\S^n\rightarrow \R^{n(n+1)/2-m}$ be the linear mapping such that the coefficient matrices $\hat{A}_1,\hat{A}_2,\ldots,\hat{A}_{n(n+1)/2-m}$ form a basis of the orthogonal complement of the span of $\{A_1,A_2,\ldots,A_m\}.$ Then the dual problem (\ref{R1tensord}) is equivalent to
\begin{equation}\label{R1tensord1}
\max\left\{ b^\top (\A\A^*)^{-1}\A(C-S):\ \hat{\A}(S)=\hat{\A}(C),\ S\in \S^n_+ \right\}.
\end{equation}
We continue to use SDPNAL+ to solve the primal problem, as it has been demonstrated in~\cite{SDPNALp1} to be highly efficient for solving the primal problem.

\begin{figure}
\centerline{\includegraphics[height=4cm,width=5cm]{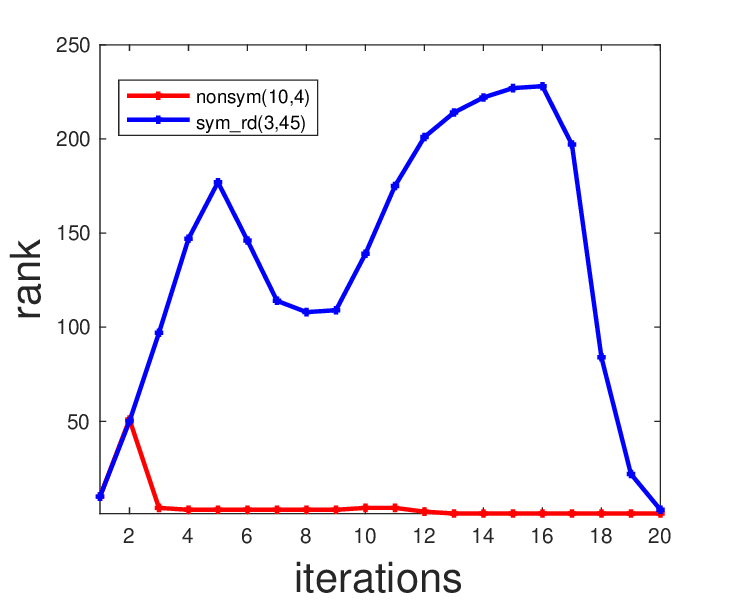}}
\caption{Rank evolution of SDPF+ for the rank-1 tensor approximation instances nonsym(10,4) and sym\_rd(3,45).} 
\label{fig:rank}
\end{figure}

\begin{center}
\begin{tiny}
\begin{longtable}{|c|ccccc|c|}
\caption{Test results for rank-1 tensor approximation SDP problems.}
\label{test-R1tensor}
\\
\hline
problem & algorithm & pfeas& dfeas & comp & fval & time \\ \hline
\endhead
nonsym(10,4)  & {\bf SDPF+}     & 9.56e-07 & 0.00e+00 & 3.03e-11 & 1.6946676e+00 & 3.26e+01 \\
n=1000,rank=999        & SDPNAL+  & 4.85e-07 & 7.29e-07 & 3.85e-05 & 1.6947107e+00 & 6.16e+01 \\
m=166374      & ManiSDP  & 3.37e-07 & 0.00e+00 & 2.72e-15 & 1.6947151e+00 & 4.70e+01 \\
\hline
nonsym(11,4)  & {\bf SDPF+}     & 2.16e-08 & 0.00e+00 & 4.42e-10 & 2.9134832e+00 & 3.41e+01 \\
n=1331,rank=1330        & SDPNAL+  & 4.66e-07 & 3.38e-07 & 3.55e-05 & 2.9134838e+00 & 1.14e+02 \\
m=287495      & ManiSDP  & 2.55e-07 & 0.00e+00 & 1.74e-15 & 2.9134851e+00 & 5.34e+01 \\
\hline
nonsym(12,4)  & {\bf SDPF+}     & 2.87e-07 & 0.00e+00 & 3.12e-11 & 5.9216048e+00 & 6.29e+01 \\
n=1728,rank=1727        & SDPNAL+  & 2.34e-07 & 0.00e+00 & 2.84e-06 & 5.9216162e+00 & 2.03e+02 \\
m=474551      & ManiSDP  & 4.28e-07 & 0.00e+00 & 4.52e-15 & 5.9216201e+00 & 1.53e+02 \\
\hline
nonsym(13,4)  & {\bf SDPF+}     & 1.92e-07 & 0.00e+00 & 3.54e-10 & 7.2745074e+00 & 1.51e+02 \\
n=2197,rank=2196        & SDPNAL+  & 8.59e-07 & 2.75e-08 & 9.28e-06 & 7.2745053e+00 & 3.88e+02 \\
m=753570      & ManiSDP  & 4.60e-07 & 0.00e+00 & 2.23e-15 & 7.2745076e+00 & 2.89e+02 \\
\hline
nonsym(14,4)  & {\bf SDPF+}     & 1.81e-08 & 0.00e+00 & 1.68e-08 & 9.6880589e+00 & 2.83e+02 \\
n=2744,rank=2743        & SDPNAL+  & 6.47e-07 & 8.27e-08 & 4.16e-05 & 9.6880401e+00 & 6.72e+02 \\
m=1157624     & ManiSDP  & 3.70e-07 & 0.00e+00 & 9.10e-16 & 9.6880544e+00 & 8.16e+02 \\
\hline
nonsym(15,4)  & {\bf SDPF+}     & 2.68e-07 & 0.00e+00 & 1.32e-08 & 1.3382864e+01 & 6.14e+02 \\
n=3375,rank=3374        & SDPNAL+  & 1.28e-07 & 3.06e-08 & 2.37e-05 & 1.3382918e+01 & 1.21e+03 \\
m=1727999     & ManiSDP  & 2.58e-07 & 4.28e-17 & 3.79e-16 & 1.3382916e+01 & 1.06e+03 \\
\hline
nonsym(16,4)  & {\bf SDPF+}     & 9.75e-08 & 0.00e+00 & 1.23e-07 & 1.5151446e+01 & 1.05e+03 \\
n=4096,rank=4095        &    SDPNAL+      &     6.87e-07     &  3.71e-07        &   3.58e-04       &    1.5151349e+01            &    2.11e+03      \\
m=2515455     &     ManiSDP     &  5.58e-03       &  0.00e+00        &   5.95e-15      &  1.1588777e+01              &   3.62e+03       \\
\hline
sym\_rd(3,45) & SDPF+     & 3.52e-07 & 7.89e-08 & 1.77e-08 & 2.1407698e+00 & 1.95e+02 \\
n=1081,rank=1080        & {\bf SDPNAL+}  & 2.60e-07 & 4.17e-07 & 1.31e-05 & 2.1407744e+00 & 6.84e+01 \\
m=211875      & ManiSDP  & 8.12e-08 & 0.00e+00 & 2.55e-12 & 2.1407693e+00 & 1.55e+03 \\
\hline
sym\_rd(3,50) & SDPF+     & 6.94e-07 & 2.10e-07 & 4.74e-09 & 2.0695003e+00 & 2.56e+02 \\
n=1326,rank=1325       & {\bf SDPNAL+}  & 9.90e-07 & 5.78e-09 & 2.74e-06 & 2.0694971e+00 & 1.38e+02 \\
m=316250      & ManiSDP  & -        & -        & -        & -             & -        \\
\hline
 sym\_rd(3,55) & SDPF+ & 5.79e-07 & 1.19e-07 & 1.80e-08 & 2.1562463e+00& 4.31e+02 \\
n=1596,rank=1595& {\bf SDPNAL+} & 8.90e-07 & 3.92e-08 & 2.98e-06 & 2.1562784e+00& 3.70e+02 \\
m=455125& ManiSDP & - & - & - & -& - \\
\hline
sym\_rd(3,60)  & SDPF+ & 4.97e-07 & 0.00e+00 & 5.30e-08 & 2.3910249e+00& 6.35e+02 \\
n=1891,rank=1890 & {\bf SDPNAL+} & 6.63e-07 & 8.73e-09 & 7.81e-07 & 2.3910013e+00& 4.15e+02 \\
m=635375 & ManiSDP & - & - &- & -& -\\
\hline
nsym\_rd([35,35,35]) & SDPF+     & 4.48e-07 & 0.00e+00 & 1.77e-08 & 3.0704739e+00 & 4.20e+02 \\
n=1225,rank=1224        & {\bf SDPNAL+}  & 7.33e-07 & 2.49e-09 & 1.80e-06 & 3.0704742e+00 & 1.25e+02 \\
m=396899      & ManiSDP  & 2.75e-09 & 0.00e+00 & 1.45e-12 & 3.0704742e+00 & 1.55e+03 \\
\hline
nsym\_rd([40,40,40]) & SDPF+     & 4.56e-08 & 0.00e+00 & 7.28e-10 & 3.8787321e+00 & 1.42e+03 \\
n=1600,rank=1599        & {\bf SDPNAL+}  & 2.17e-07 & 1.43e-08 & 8.73e-06 & 3.8787395e+00 & 1.62e+02 \\
m=672399      & ManiSDP  & 2.06e-05 & 1.44e-05 & 2.28e-08 & 3.8785974e+00 & 3.60e+03 \\
\hline
\end{longtable}
\end{tiny}
\end{center}
From Table~\ref{test-R1tensor}, we observe that SDPF+ is more efficient than ManiSDP and achieves better accuracy than SDPNAL+. For the first 7 instances, SDPF+ outperforms SDPNAL+ in efficiency, as these problems have rank-1 dual solutions, and the ranks of the iterates within the ALM subproblems remain low. For instances 8–13, SDPF+ is slower than SDPNAL+ because the ranks of the iterates increase significantly during the early stage of the ALM iterations before eventually decreasing. See Figure~\ref{fig:rank} on the rank update of  SDPF+ for the instances 
nonsym(10,4) and sym\_rd(3,45). Nevertheless, it is worth noting that SDPF+ is the only solver capable of solving all instances to the required accuracy.

\subsection{Second order reduced density matrices}
We test SDP problems arising from fermion second-order reduced density matrices (2RDM)~\cite{nakata2008variational, nakata2001variational} for
molecules. 
For each molecule, there are four types of relaxations available\footnote{The dataset is available at {\tt https://nakatamaho.riken.jp/rdmsdp/sdp\_rdm.html}.}, and we select the tightest one -- whose name ends with ``1t2p'' -- for testing. This set of problems is particularly challenging due to the following reasons:

\begin{itemize}
\item They are multi-block SDP problems, with more than 20 blocks in most instances.
\item The coefficient matrices in the constraints are significantly denser than those in previous problems, and the matrix $M_R$ defined in~(\ref{linsolveM}) is nearly fully dense.
\end{itemize}

Although the number of constraints in these problems is not as large as in the rank-1 tensor approximation problems, we do not compare with SDPT3 because it requires more than 10 minutes per iteration, largely due to the density of the coefficient matrices. Similarly, we exclude ManiSDP from the comparison
because based on our experiments, it is unable to solve these problems to the required accuracy.

We set the KKT residual tolerance to $10^{-5}$, since neither SDPF+ nor SDPNAL+ can solve many of the instances to $10^{-6}$ within 3600 seconds. Additionally, we only test on instances with fewer than 16,000 constraints to ensure that both SDPF+ and SDPNAL+ can successfully solve most of them within the given time limit.

\begin{figure}
\centerline{\includegraphics[height=4cm,width=12cm]{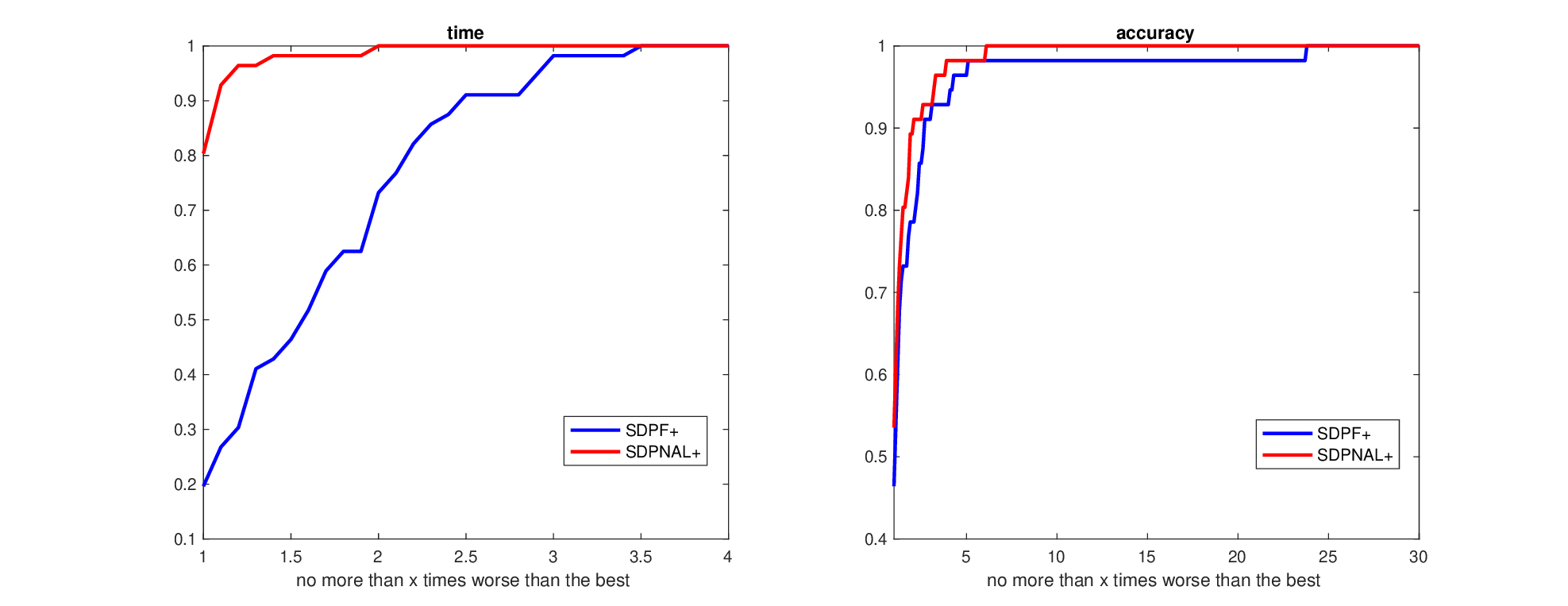}}
\caption{Performance profile for the test results of 2RDM.} 
\label{fig:2RDM}
\end{figure}

From Figure~\ref{fig:2RDM} and Table~\ref{test-2RDM} in Appendix A, we observe that SDPF+ is slightly slower than SDPNAL+, with both solvers achieving a similar level of accuracy. One contributing factor is that the matrix $M_R$ is dense, requiring SDPF+ to use dense Cholesky decomposition to compute the preconditioner. This step is computationally expensive and produces a dense Cholesky factor, which slows down the forward and backward substitution operations. Nevertheless, the performance gap is within an acceptable range, as the run time of SDPF+ remains within four times that of SDPNAL+ for most instances.

\subsection{Hans Mittelmann benchmark dataset}

We test the Hans Mittelmann benchmark suite for SDP problems\footnote{\href{https://plato.asu.edu/ftp/sdp/}{https://plato.asu.edu/ftp/sdp/}}. This dataset contains a wide variety of linear SDP problems with diverse structures. Since some problems are small scale but extremely ill-conditioned and challenging to solve, they are 
usually only solvable by interior-point methods. 

We do not compare with ManiSDP, as it is highly sensitive to its parameter settings. For example, for the Lov\'asz Theta SDP problems, we had to carefully tune ManiSDP’s parameters to achieve convergence; using default parameters, the solver rarely converges.

For some instances, the duality gap is much smaller than the complementarity, so we report $\min\{\text{pdgap}, \text{comp}\}$ in the Table~\ref{test-HMT} for clarity. To avoid redundant testing, we exclude the Lov\'asz Theta SDP problems 
that are already tested earlier in Section 4.1. 

The dataset also contains sparse SDP problems with hundreds or even thousands of blocks, each with a very small block size, for example, the {\tt mater-*} instances. These problems are particularly well-suited for interior-point methods, and SDPT3 significantly outperforms both SDPF+ and SDPNAL+ on such problems. Since SDPF+ and SDPNAL+ are designed for solving SDPs with large matrix dimensions rather than a large number of small blocks, we omit testing on this type of problems. In total, we tested 71 instances.

\begin{figure}[H]
\centerline{\includegraphics[height=4cm,width=12cm]{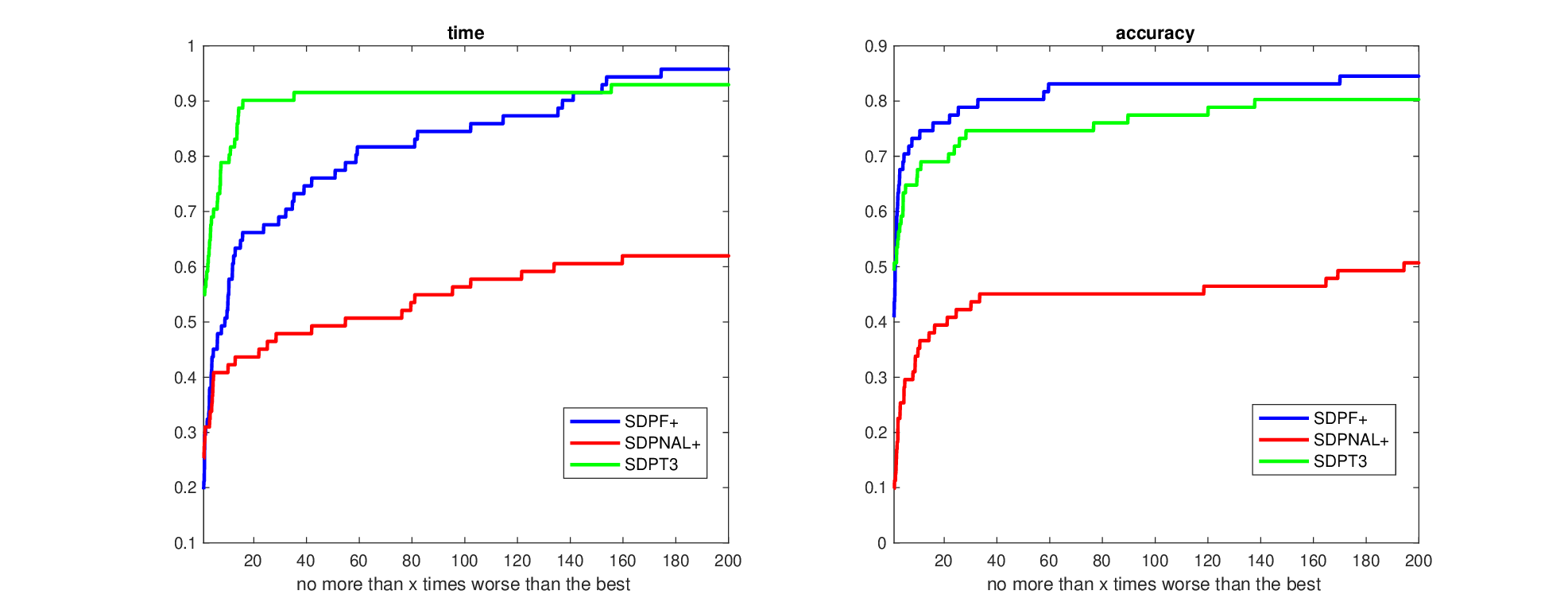}}
\caption{Performance profile for test results on Hans Mittelmann benchmark dataset.} 
\label{fig:HMT}
\end{figure}

\begin{center}
\begin{small}
\begin{longtable}{|lccc|}
\caption{Percentage of problems solved by each solver to different accuracy levels}
\\
\hline
\label{tab:solver-accuracy}
Solver & $\leq 10^{-6}$ & $\leq 10^{-4}$ & $\leq 10^{-2}$ \\
\hline
SDPF+    & 71.83\% & 83.10\% & 95.77\% \\
SDPNAL+  & 23.94\% & 54.93\% & 66.20\% \\
SDPT3    & 56.34\% & 80.28\% & 91.55\% \\
\hline
\end{longtable}
\end{small}
\end{center}

From Figure~\ref{fig:HMT} and Table~\ref{test-HMT} in Appendix A, we observe that SDPT3 is the most efficient solver, with SDPF+ performing in between SDPT3 and SDPNAL+ in terms of efficiency. This result aligns with our expectations, as interior-point methods are known to be efficient and robust for small but highly ill-conditioned problems -- characteristics that are common in many instances from the Hans Mittelmann benchmark suite.

Interestingly, as shown in Table~\ref{tab:solver-accuracy}, SDPF+ solves the largest number of instances to the accuracy level of $10^{-6}$, which exceeds our expectations. One possible explanation is that SDPT3 tends to terminate prematurely when there is little progress or when the solver runs into numerical difficulties due to extreme ill-conditioning. Moreover, for the {\tt G**\_mb} and {\tt G**mc} instances, SDPF+ can
be more than 5--10 times faster than SDPT3, and even more than 100 times
faster than SDPNAL+.
However, for some difficult instances such as 
{\tt trto*}, {\tt vibra*}, SDPF+ can be more than 10 times slower than 
SDPT3, while SDPNAL+ fails to solve them.

\subsection{Non-parametric density estimation}

Here we test nonlinear SDP problems arising from
non-parametric models of density estimation \cite{marteau2020non}: 
\begin{equation}\label{DESDP}
\min\left\{ -\frac{{\bf 1}_{1\times n}}{n} \log \( \dd(VXV^\top) \)+\lambda_1\<I_n,X\>+\lambda_2 \|X\|^2_F:\ \<M,X\>=1,\ X\in \S^n_+ \right\},
\end{equation}
where $V=KL^{-1}$ with $K=\{\kappa(x_i,x_j)\}_{n\times n}$ for some kernel function $\kappa$ and sample points $\{x_i\}_{i\in [n]}$ from the distribution that we want to estimate; $L$ is the lower Cholesky factor of $K$; $M\in \S^n_+$ is a matrix that is used to normalize the distribution (see \cite{marteau2020non} for its formula). 

We compare SDPF+ with FISTA~\cite{beck2009fast}, which was used in~\cite{marteau2020non} to solve the dual problem of~(\ref{DESDP}). Following the numerical setup in~\cite{marteau2020non}, we randomly sample $n$ points from a 5-dimensional mixed Gaussian distribution:
\begin{equation}\label{gaussiandis}
0.5 \cdot \mathcal{N}(\mu_1, I_d / (2\pi)) + 0.5 \cdot \mathcal{N}(\mu_2, I_d / (2\pi)),
\end{equation}
where $\mu_1 = -2e_1$ and $\mu_2 = 2e_1$, with $e_1 \in \mathbb{R}^d$ being the standard basis vector whose first entry is $1$ and the rest are $0$. We use the Gaussian kernel $\kappa(x, y) := \exp\left(-\|x - y\|^2 / (2\sigma^2)\right)$.

In~\cite{marteau2020non}, the parameters $\lambda_1$, $\lambda_2$, and $\sigma$ are selected via cross-validation to optimize model accuracy. Since our goal is to evaluate algorithmic efficiency rather than prediction accuracy, we bypass cross-validation and directly set $\lambda_1 = 0.1$ and $\lambda_2 = 0.005$, which typically yield good performance. As for $\sigma$, the optimal choice is $\sigma = 1/\sqrt{2\pi}$, which aligns with the covariance of the mixed Gaussian distribution.

We observe that small deviations from the optimal $\sigma$, such as $\sigma = 2/\sqrt{2\pi}$ or $\sigma = 0.5/\sqrt{2\pi}$, can significantly impact both the prediction accuracy and the difficulty of the resulting optimization problem. Therefore, we consider $\sigma \in \{0.5/\sqrt{2\pi}, 1/\sqrt{2\pi}, 2/\sqrt{2\pi}\}$ and $n \in \{500, 1000, 2000\}$ to evaluate the scalability and robustness of the algorithms.

\begin{center}
\begin{tiny}
\begin{longtable}{|c|ccccc|c|}
\caption{Test results for non-parametric density estimation.}
\\
\hline
\label{test-DE}
problem & algorithm & pfeas& dfeas & comp & fval & time \\ \hline
\endhead
n=500, rank=10 & {\bf SDPF+} & 6.43e-09 & 1.32e-07 & 5.79e-09 & 1.0215680e+01& 2.03e+01 \\
$\sigma = 0.5/\sqrt{2\pi}$ & FISTA & 9.91e-07 & 2.36e-08 & 3.68e-08 & 1.0215698e+01& 3.36e+01 \\
 \hline
n=500, rank=5 & {\bf SDPF+} & 1.83e-09 & 9.62e-07 & 1.31e-08 & 3.1095986e+00& 1.07e+01 \\
$\sigma = 1/\sqrt{2\pi}$ & FISTA & 2.36e-07 & 7.46e-07 & 3.38e-07 & 3.1095982e+00& 4.65e+01 \\
\hline
 n=500, rank=39 & {\bf SDPF+} & 6.49e-09 & 1.75e-07 & 2.26e-07 & 3.2469466e+00& 1.06e+01 \\
$\sigma = 2/\sqrt{2\pi}$ & FISTA & 1.54e-07 & 5.33e-07 & 3.52e-07 & 3.2469462e+00& 3.60e+03 \\
\hline
  n=1000, rank=11 & {\bf SDPF+} & 2.17e-09 & 6.68e-07 & 3.33e-08 & 9.1833437e+00& 3.19e+01 \\
$\sigma = 0.5/\sqrt{2\pi}$ & FISTA & 9.59e-07 & 4.84e-08 & 5.77e-08 & 9.1833572e+00& 1.46e+02 \\
\hline
n=1000, rank=8& {\bf SDPF+} & 3.10e-10 & 2.01e-07 & 3.48e-08 & 3.1707847e+00& 3.77e+01 \\
$\sigma = 1/\sqrt{2\pi}$ & FISTA & 1.07e-07 & 0.00e+00 & 2.48e-07 & 3.1707849e+00& 3.47e+02 \\
 \hline
n= 1000, rank=38 & {\bf SDPF+} & 5.50e-14 & 4.44e-07 & 2.08e-07 & 3.2405354e+00& 3.75e+01 \\
$\sigma = 2/\sqrt{2\pi}$ & FISTA & - & - & - & -& - \\
 \hline
  n=2000, rank=8 & {\bf SDPF+} & 1.60e-09 & 3.37e-07 & 2.95e-09 & 8.4964317e+00& 1.28e+02 \\
$\sigma = 0.5/\sqrt{2\pi}$ & FISTA & 9.61e-07 & 6.17e-08 & 8.08e-08 & 8.4964428e+00& 7.36e+02 \\
\hline
n=2000, rank=7 & {\bf SDPF+} & 4.92e-09 & 6.69e-07 & 9.01e-10 & 3.2101141e+00& 1.83e+02 \\
$\sigma = 1/\sqrt{2\pi}$ & FISTA & 8.58e-08 & 0.00e+00 & 2.51e-07 & 3.2101142e+00& 2.72e+03 \\
\hline
n=2000, rank=30 & {\bf SDPF+} & 4.44e-16 & 7.77e-07 & 4.58e-08 & 3.2109010e+00& 3.87e+02 \\
$\sigma = 2/\sqrt{2\pi}$ & FISTA & - & - & - & -& - \\
\hline
\end{longtable}
\end{tiny}
\end{center}

From Table~\ref{test-DE}, we observe that SDPF+ is more efficient than FISTA across all tested instances. In particular, for the third instance, SDPF+ is over 100 times faster than FISTA. This performance advantage is primarily due to the low-rank structure of the problems, which aligns well with the design of the low-rank solver SDPF+. Moreover, SDPF+ exhibits a more stable performance than FISTA with respect to variations in the kernel bandwidth parameter $\sigma$.

\subsection{Weighted nearest correlation matrix}

In this subsection, we test on quadratic SDP problems arising from H-weighted nearest correlation matrix estimation problem \cite{qi2011augmented}:
\begin{equation}\label{HNC}
\min\left\{ 
\frac{1}{2}\|H\circ (X-\hat{X})\|_F^2 
\ :\ 
\dd(X) = \mathbf{1}_n,\ 
X \in \mathbb{S}_+^n 
\right\},
\end{equation}
where $\hat{X} \in \mathbb{S}^n$ is the given data matrix. 
This formulation corresponds to a weighted projection of $\hat{X}$ 
onto the set of correlation matrices, where the weight matrix $H$ 
assigns entrywise importance in the Frobenius norm.
We compare SDPF+ with the dual augmented Lagrangian method with subproblem solved by a semismooth Newton method (SSNAL), developed by Qi and Sun  in~\cite{qi2011augmented}. The {\sc Matlab}
implementation of SSNAL is available online\footnote{\href{https://www.polyu.edu.hk/ama/profile/dfsun//Codes/correlation-matrix/}{https://www.polyu.edu.hk/ama/profile/dfsun//Codes/correlation-matrix/}}. We use the data generation procedure provided in the  code to generate the matrices $H$ and $\hat{X}$. The problem dimension $n$ is chosen from the set $\{100, 200, 500, 1000, 2000, 5000\}$.

\begin{center}
\begin{tiny}
\begin{longtable}{|c|ccccc|c|}
\caption{Test results for H-weighted nearest correlation matrix problem with square loss function.}
\\
\hline
problem & algorithm & pfeas& dfeas & comp & fval & time \\ \hline
\endhead
\label{test-HNC}
n = 100  & {\bf SDPF+} & 9.94e-09 & 8.63e-07 & 2.33e-08 & 7.6275514e+01& 8.63e+00 \\
rank=56 & SSNAL & 8.68e-10 & 2.00e-08 & 6.16e-08 & 7.6275517e+01& 1.27e+01 \\
\hline
  n = 200 & SDPF+ & 9.26e-09 & 4.01e-07 & 1.73e-08 & 4.4762000e+02& 1.73e+01 \\
rank= 100 & {\bf SSNAL} & 6.58e-10 & 1.05e-08 & 2.58e-08 & 4.4762000e+02& 8.73e+00 \\
\hline
n = 500   & SDPF+ & 4.44e-09 & 8.36e-07 & 2.71e-09 & 3.7958633e+03& 3.40e+01 \\
rank=222 & {\bf SSNAL} & 1.39e-10 & 1.96e-09 & 9.51e-09 & 3.7958627e+03& 2.40e+01 \\
 \hline
 n = 1000  & SDPF+ & 8.32e-09 & 3.81e-07 & 4.60e-10 & 1.8637052e+04& 7.73e+01 \\
rank=396 & {\bf SSNAL} & 2.51e-10 & 4.88e-10 & 9.51e-09 & 1.8637051e+04& 6.92e+01 \\
 \hline
 n = 2000  & {\bf SDPF+} & 9.13e-09 & 3.32e-07 & 8.04e-12 & 9.0980052e+04& 3.06e+02 \\
rank=697& SSNAL & 2.44e-11 & 1.95e-10 & 9.36e-09 & 9.0980052e+04& 3.31e+02 \\
\hline
  n = 5000 & {\bf SDPF+} & 1.47e-11 & 6.96e-07 & 2.41e-08 & 7.1678529e+05& 2.36e+03 \\
 rank=1441& SSNAL & 1.00e-10 & 5.50e-10 & 2.51e-08 & 7.1678527e+05& 3.26e+03 \\
\hline
\end{longtable}
\end{tiny}
\end{center}

From Table~\ref{test-HNC}, we observe that the run of SDPF+ is roughly comparable to that of SSNAL. Since SSNAL is a specialized solver tailored for problem~(\ref{HNC}), and the solution ranks of the tested instances are rather high, these results sufficiently demonstrate the efficiency of SDPF+ in handling problem~(\ref{HNC}).

Next, we change the problem (\ref{HNC}) by replacing the square-loss function by the Huber-loss function \cite{huber1992robust}:
\begin{equation}\label{HNCHB}
\min\left\{f\(H\circ (X-\hat{X})\):\ \dd(X)={\bf 1}_n,\ X\in \S^n_+\right\},
\end{equation}
where $f(X):=\sum_{i,j\in [n]} g(X_{ij})$ and $g:\R\rightarrow \R$ is defined as
\begin{equation}\label{defihuber}
g(x):=\begin{cases} \frac{x^2}{2} & |x|\leq \delta \\ \delta\cdot \(x-\frac{\delta}{2}\) & |x|> \delta. \end{cases}
\end{equation}
Since SSNAL cannot handle the Huber loss function, we compare SDPF+ with the accelerated projected gradient (APG) method~\cite{jiang2012inexact}. We use the {\sc Matlab} function \texttt{nearcorr} to compute the projection onto the feasible region of problem~(\ref{HNCHB}) with an accuracy tolerance of $10^{-8}$. We set $\delta = 0.1$ and generate the data matrices in the same manner as in the previous experiment.

\begin{center}
\begin{tiny}
\begin{longtable}{|c|ccccc|c|}
\caption{Test results for H-weighted nearest correlation matrix problems with the Huber-loss function.}
\\
\hline
\label{test-HNCHB}
problem & algorithm & pfeas& dfeas & comp & fval & time \\ \hline
\endhead
 n=100 & {\bf SDPF+} & 8.16e-09 & 7.22e-07 & 1.74e-07 & 5.0989963e+01& 2.21e+01 \\
 rank=56& APG & 0.00e+00 & 9.79e-07 & 2.68e-16 & 5.0989960e+01& 6.44e+01 \\
\hline
 n=200 & {\bf SDPF+} & 9.42e-09 & 7.42e-07 & 8.07e-08 & 2.6947628e+02& 2.65e+01 \\
rank=99& APG & 0.00e+00 & 9.90e-07 & 3.50e-16 & 2.6947628e+02& 1.85e+02 \\
\hline
n=500 & {\bf SDPF+} & 1.19e-11 & 7.82e-07 & 3.47e-08 & 2.1158843e+03& 7.69e+01 \\
 rank=216 & APG & 0.00e+00 & 9.61e-07 & 2.97e-16 & 2.1158842e+03& 1.35e+03 \\
 \hline
  n=1000& {\bf SDPF+} & 4.44e-09 & 9.46e-07 & 8.16e-08 & 9.8461211e+03& 1.84e+02 \\
 rank=384 & APG & 0.00e+00 & 2.32e-04 & 1.39e-15 & 9.8461306e+03& 3.60e+03 \\
 \hline
n=2000  & {\bf SDPF+} & 5.29e-10 & 6.54e-07 & 7.68e-08 & 4.5499748e+04& 8.53e+02 \\
rank=654& APG & 0.00e+00 & 3.75e-03 & 6.59e-17 & 4.5510046e+04& 3.60e+03 \\
\hline
n=5000  & {\bf SDPF+} & 3.47e-09 & 8.27e-04 & 9.82e-03 & 3.3638798e+05& 3.61e+03 \\
rank=1465 & APG & - & - & - & -& - \\
\hline
\end{longtable}
\end{tiny}
\end{center}

From Table~\ref{test-HNCHB}, we observe that SDPF+ is more efficient than APG across all tested instances. This result highlights both the efficiency and generality of SDPF+.

\section{Conclusion}\label{Sec-conc}

In this work, we propose a preconditioned ALM that employs a weighted penalty function as the preconditioner,  
where the weight matrix is derived from the projection operator onto the tangent space of the feasible set.  
Building on this algorithm, we further develop SDPF+, a two-stage SDP solver that combines a primal feasible method with the PALM.  
The solver is designed to handle nonlinear objective functions and general linear convex constraints.  
Extensive numerical experiments on a broad range of benchmark problems demonstrate the robustness, efficiency, and versatility of SDPF+.  
The results show that SDPF+ performs competitively against several state-of-the-art solvers, validating its effectiveness as a general-purpose SDP solver.  
As a future direction, it would be valuable to explore more memory-efficient preconditioning techniques for large-scale problems and to extend the preconditioned framework for handling inequality and other conic constraints more effectively. It is worth noting that SDPF+ is generally outperformed by SDPNAL+ in our numerical experiments when the optimal solutions have high rank. This behavior is expected, as SDPF+ relies on low-rank factorization, and its computational advantage diminishes when the iterates have high rank. An interesting direction for future research is whether similar preconditioning strategies can be incorporated to accelerate traditional convex SDP algorithms such as ADMM and augmented Lagrangian methods.



\bibliographystyle{abbrv}
\bibliography{SDPF+}

\newpage
\appendix

\begin{small}
\section{Proof of Theorem~\ref{convtheo}}
\label{Proof-of-theorem-1}

In the following analysis, we always assume $k\geq 2$. From (\ref{defig}) and (\ref{defiH}), we have that
\begin{equation}\label{deveg}
 \partial_{z} g(z^k)-\H^*(\gamma^k)= \(\partial_X \delta_{\S^n_+}(X^k)+\nabla f(X^k)-\A^*(\lambda^k)-\B^*(\mu^k), \partial_y \delta_\P(y^k)+\mu^k \).
\end{equation}
Because $X^k=R^kR^{k\top},$ $X^k\in \S^n_+.$ From (\ref{ALMsubYZcl}) and the dual update of $\mu$ in Algorithm~\ref{alg:PALMG}, we know that
\begin{equation}\label{musatis}
0\in \partial_y \delta_\P(y^k)+\mu^k .
\end{equation}
From Assumption~\ref{assALMsub} and the dual update in Algorithm~\ref{alg:PALMG}, we have that 
\begin{equation}\label{lamsatis}
\dist\( \partial_X \delta_{\S^n_+}(X^k)+\nabla f(X^k)-\A^*(\lambda^k)-\B^*(\mu^k),0\)<\epsilon_{k-1}. 
\end{equation}
Combining (\ref{musatis}) and (\ref{lamsatis}), we have that $\dist\( \partial_{z} g(z^k)-\H^*(\gamma^k),0 \)\rightarrow 0.$ Thus, we only have to prove that $\|\H(z^k)-b\|_2\rightarrow 0.$ The equations (\ref{musatis}) and (\ref{lamsatis}) imply that there exist $v^k\in \X$ such that 
\begin{equation}\label{condv}
v^k\in \partial_{z} g(z^k)-\H^*(\gamma^k),\ \|v^k\|_2< \epsilon_{k-1}. 
\end{equation}
From (\ref{condv}) and the convexity of $g,$ we have that
\begin{equation}\label{desine}
g(z^{k+1})\geq g(z^k)+\< v^k+\H^*(\gamma^k),z^{k+1}-z^k \>. 
\end{equation}
Define 
\begin{equation}\label{defiL}
L_k:= g(z^k)-\<\gamma^k,\H(z^k)-c\>.
\end{equation}
From (\ref{desine}), we have that  
\vspace{-5mm}
\begin{align}
L_{k+1}-L_k &\geq -\<\gamma^{k+1},\H(z^{k+1})-c\>+\<\gamma^{k},\H(z^{k})-c\>+\< v^k+\H^*(\gamma^k),z^{k+1}-z^k \> 
\notag \\
&=\< \gamma^k-\gamma^{k+1},\H(z^{k+1})-c \>+\<v^k,z^{k+1}-z^k\>.
\label{Aug_14_1}
\end{align}
From the dual update of Algorithm~\ref{alg:PALMG}, we have that
\begin{align}
&\< \gamma^k-\gamma^{k+1},\H(z^{k+1})-c \>
=\<\lambda^k-\lambda^{k+1},\A(X^{k+1})-b\>+\<\mu^k-\mu^{k+1},\B(X^{k+1})-y^{k+1}\> \notag 
\\
&=\beta_1\| \A(X^{k+1})-b \|_{W^k}^2+\beta_2\| \B(X^{k+1})-y^{k+1} \|_2^2\notag 
\\
&\geq \beta_1 \sigma\| \A(X^{k+1})-b \|_2^2+\beta_2\| \B(X^{k+1})-y^{k+1} \|_2^2, \label{Aug_14_2}
\end{align}
where the last inequality follows from that $W^k\succeq \sigma I_m.$ Substituting (\ref{Aug_14_2}) into (\ref{Aug_14_1}), we get
\begin{equation}\label{Aug_14_2a}
L_{k+1}\geq L_k+\beta_1 \sigma\| \A(X^{k+1})-b \|_2^2+\beta_2\| \B(X^{k+1})-y^{k+1} \|_2^2+\<v^k,z^{k+1}-z^k\>. 
\end{equation}
By applying (\ref{Aug_14_2a}) from $k=2$ to $K\in \N^+,$ we get 
\begin{equation}\label{Aug_14_3}
L_{K+1}\geq L_2+\sum_{k=2}^K\(\beta_1 \sigma\| \A(X^{k+1})-b \|_2^2+\beta_2\| \B(X^{k+1})-y^{k+1} \|_2^2+\<v^k,z^{k+1}-z^k\>\). 
\end{equation}
From  Assumption~\ref{assFeas}, we have that 
\begin{equation}\label{Aug_14_4}
g(z^*)=g(z^*)-\<\gamma^{K+1},\H(z^*)-c\>.
\end{equation}
From (\ref{condv}), we have that
\begin{equation}\label{Aug_14_5}
g(z^*)-\<\gamma^{K+1},\H(z^*)-c\>\geq g(z^{K+1})-\<\gamma^{K+1},\H(z^{K+1})-c\>+\<v^{K+1},z^*-z^{K+1}\>.
\end{equation}
Using the definition~(\ref{defiL}) and (\ref{Aug_14_4}), (\ref{Aug_14_5}) can be simplified as
\begin{equation}\label{Aug_14_6}
g(z^*)\geq L_{K+1}+\<v^{K+1},z^*-z^{K+1}\>.
\end{equation}
Substituting (\ref{Aug_14_3}) into (\ref{Aug_14_6}), we get
\begin{align}\label{Aug_14_7}
&g(z^*)-\<v^{K+1},z^*-z^{K+1}\>\notag \\
&\geq L_2+\sum_{k=2}^K\(\beta_1 \sigma\| \A(X^{k+1})-b \|_2^2+\beta_2\| \B(X^{k+1})-y^{k+1} \|_2^2+\<v^k,z^{k+1}-z^k\>\). 
\end{align}
Because $\{z^k\}$ is bounded, there exists $M>0$ such that $\|z^k\|_2<M$ for any $k\in \N^+.$ Using this as well as $\|v^k\|_2\leq \epsilon_{k-1}$ in (\ref{Aug_14_7}), we get
\vspace{-5mm}
\begin{align}\label{Aug_14_8}
&g(z^*)-L_2+\epsilon_K\(\|z^*\|+M\)+2M\sum_{k=1}^{K-1}\epsilon_k\notag \\
&\geq \sum_{k=2}^K\(\beta_1 \sigma\| \A(X^{k+1})-b \|_2^2+\beta_2\| \B(X^{k+1})-y^{k+1} \|_2^2\). 
\end{align}
From the summability of $\{\epsilon_k\}$ in Assumption~\ref{assALMsub}, the left hand side of (\ref{Aug_14_8}) is bounded. Therefore, we have that $\|\H(z^k)-c\|_2\rightarrow 0.$ This completes the proof of Theorem~\ref{convtheo}.

\section{Usage and data structures of SDPF+}\label{App-solver}
In this Section, we describe the usage and data structures of SDPF+.  
Although the constraint $\H(X) \in \P$ in (\ref{SDPG}) encompasses many types of convex constraints, the implementation of SDPF+ includes several special cases that can be specified by users more conveniently.

Our solver can handle the following multi-block version of \eqref{SDPG}:
\begin{equation}\label{SDPG-multi}
\begin{array}{rl}
\min & f(X) \\[3pt]
\mbox{s.t.} & \A(X)=b,\ l\leq \B(X)\leq u,\\[3pt]
& X\in \P,\ \H(X)\in \Q,\quad
 X=(X_1,\ldots,X_K)\in \cV^{(1)}\times\cdots\times \cV^{(K)},
\end{array}
\tag{SDPG1}
\end{equation}
where for each $i\in [K]$, $\cV^{(i)}\subset \V^{(i)}$ is either $\S^{n_i}_+\subset \S^{n_i}$, $\R^{n_i}_+\subset \R^{n_i}$, or 
$\R^{n_i}\subset \R^{n_i}$. The mappings $\A(X),$ $\B(X)$ and $\H(X)$ are given by 
\vspace{-5mm}
\begin{align}
 &\A(X) = \A^{(1)}(X_1) + \cdots + \A^{(K)}(X_K), \quad
   \B(X) = \B^{(1)}(X_1) + \cdots + \B^{(K)}(X_K), \notag \\
 & \H(X) = \H^{(1)}(X_1) + \cdots + \H^{(K)}(X_K), \notag
\end{align}
where for each $i\in [K]$, $\A^{(i)}:\V^{(i)}\to \R^m,$ $\B^{(i)}: \V^{(i)}\to \R^q$ 
and $\H^{(i)}:\V^{(i)}\to \R^p$ are given linear mappings. Moreover, $\P$ is a closed convex subset of $\V^{(1)}\times \cdots\times \V^{(K)}$, $\Q$ is a closed convex subset of $\R^p$, $b\in \R^m$, $l,u$ are $q$-dimensional vectors whose entries could be $-\infty$ and $\infty$, respectively.

The main routine of SDPF+ is \texttt{sdpfplus.m} with the following calling syntax:
\begin{center}
\texttt{runhist = sdpfplus(blk,At,C,b,Bt,l,u,par);}
\end{center}
The input arguments of SDPF+ are described as follows:
\begin{enumerate}
\item \texttt{blk} is a $K \times 2$ cell array that specifies the block structure of the variables, where $K$ is the number of blocks. For each $k \in [K]$, \texttt{blk\{k,1\}} is a character indicating the type of the $k$-th block, which can be one of the following:
\begin{itemize}
\item \texttt{'s'}: symmetric positive semidefinite matrix variable,
\item \texttt{'l'}: nonnegative vector variable,
\item \texttt{'f'}: free vector variable.
\end{itemize}
\texttt{blk\{k,2\}} specifies the corresponding dimension of the $k$-th block.

\item \texttt{At} is a $K\times 1$ or $K\times 2$ cell array that describe the linear operator $\A$ in (\ref{SDPG-multi}). It has two columns when the user wants to make use of the sparse plus low-rank structure of the constraint matrices $A_i$'s described in  \cite{tang2024exploring}. Suppose the $k$-th block is a matrix variable with dimension $n$ and the constraint matrices 
$\{A_i\}_{i\in[m]}$ have sparse plus low-rank decomposition of the form: $A_i=A_i^s+U_i\dd(d_i)U_i^\top$ such that $d_i\in \R^{r_i}$ for any $i\in [m].$ Then \texttt{At\{k,1\}} is a $n(n+1)/2\times m$ matrix such that for any $i\in [m]$ and $1\leq \alpha\leq \beta\leq n:$
\begin{equation}\label{At1}
\texttt{At\{k,1\}}[(\beta-1)\beta/2+\alpha,i]=\begin{cases} A_i^s[\alpha,\beta] & {\rm if}\ \alpha = \beta \\[3pt] 
\sqrt{2}A_i^s[\alpha,\beta] & {\rm if}\ \alpha<\beta,
\end{cases}
\end{equation}
\texttt{At\{k,2\}} is a $(n+2)\times \sum_{i=1}^mr_i$ matrix such that:
\begin{equation}\label{At2}
\texttt{At\{k,2\}}=\begin{bmatrix}
\begin{bmatrix}
1\cdot{\bf 1}_{1\times r_1}\\
U_1 \\
d_1^\top
\end{bmatrix}
&
\begin{bmatrix}
2\cdot{\bf 1}_{1\times r_2}\\
U_2 \\
d_2^\top
\end{bmatrix}
&
\cdots
\begin{bmatrix}
m\cdot{\bf 1}_{1\times r_m}\\
U_m \\
d_m^\top
\end{bmatrix}
\end{bmatrix}.
\end{equation}
When the low-rank structure of the constraint matrices is not used, $\texttt{At}$ will only have one column.

\item \texttt{C} is either a cell array or a function handle that describes the objective function $f$ in (\ref{SDPG-multi}). When the SDP problem has linear objective function $f(X)=\sum_{i=1}^K\<C_i,X_i\>$ such that the matrix blocks have sparse plus low-rank decomposition $C_i=C_i^s+U_i\dd(d_i)U_i^\top,$ we have that
\begin{equation}\label{C1}
\texttt{C}\{i,1\}=\begin{cases}
C_i & {\rm if}\ \texttt{blk\{i,1\}}\neq \texttt{'s'} \\
C_i^s & {\rm if}\ \texttt{blk\{i,1\}}= \texttt{'s'}.
\end{cases}
\end{equation}
\begin{equation}\label{C2}
\texttt{C}\{i,2\}=\begin{cases}
\texttt{[]} & {\rm if}\ \texttt{blk\{i,1\}}\neq \texttt{'s'} \\
\begin{bmatrix} U_i\\d_i^\top \end{bmatrix} & {\rm if}\ 
\texttt{blk\{i,1\}}= \texttt{'s'}.
\end{cases}
\end{equation}
When $f(X)$ is a nonlinear function, SDPF+ input $f(\cdot)$ using a function handle. The input to the function handle is a cell array representing the primal variable, i.e., $\{X_1; \ldots; X_K\}$, and the output consists of the function value $f(X)$ and a cell array containing the gradients, i.e., $\{\nabla_{X_1} f(X); \ldots; \nabla_{X_K} f(X)\}$. Consider the following example:
\begin{equation}\label{exfun}
f(X)=\<I_n,X_1\>^2+\gamma\exp\({\bf 1}_\ell^\top X_2\),
\end{equation}
where $X_1\in \S^n$, $X_2\in \R^{\ell}$ and $\gamma>0$ is a parameter. We can define $\texttt{C = @Cfun(X,gama)},$ where $\texttt{Cfun}$ is defined 
in Listing~\ref{lst:nonobj}.

In some cases, the function value $f(RR^\top)$ and its gradient $\nabla f(RR^\top)$ can be computed more efficiently without explicitly forming the full matrix $RR^\top$ for each matrix block. To take advantage of this, users may define the input $X$ as a \texttt{struct}-type array with a field named \texttt{R} representing the low-rank factor. Then \texttt{R} is a $K\times 1$ cell array such that for any $i\in [K]$
\begin{equation}\label{XR}
\texttt{R\{i\}}=\begin{cases}
U\ {\rm such\ that}\ UU^\top=X_i & {\rm if}\ \texttt{blk\{i,1\}}=\texttt{`s'} \\
X_i & {\rm if}\ \texttt{blk\{i,1\}}\neq \texttt{`s'}. 
\end{cases}
\end{equation}
We again use the example (\ref{exfun}) to present a  computationally more
efficient way to define the function in Listing~\ref{lst:nonobj1}.

\medskip
\begin{lstlisting}[caption={Example for defining a nonlinear objective function.}, label={lst:nonobj}] 
function [f,g] = Cfun(X,gama)
n = size(X{1},1);
l = size(X{2},1);
f = (trace(X{1}))^2+gama*exp(sum(X{2}));
g = cell(2,1);
g{1} = 2*trace(X{1})*speye(n);
g{2} = gama*ones(l,1)*exp(sum(X{2}));
end
\end{lstlisting}
\begin{lstlisting}[caption={Example for defining a nonlinear objective function using the low-rank factor.}, label={lst:nonobj1}]
function [f,g] = Cfun(X,gama)
R = X.R;
n = size(R{1},1);
l = size(R{2},1);
f = R(:)'*R(:)+gama*exp(sum(R{2}));
g = cell(2,1);
g{1} = 2*R(:)'*R(:)*speye(n);
g{2} = gama*ones(l,1)*exp(sum(R{2}));
end
\end{lstlisting}

\item \texttt{b} is an $m$-dimensional vector which is exactly the vector $b$ in (\ref{SDPG-multi}). 

\item \texttt{Bt,l,u} are used to denote the inequality constraint, $l\leq \B(X)\leq u$, in (\ref{SDPG-multi}). The data structure of \texttt{Bt} is the same as \texttt{At}.

\item \texttt{par} is a \texttt{struct}-type array that contains the parameters and options for the solver. The user can refer to the main function \texttt{sdpfplus.m} for the complete list of available parameters. Although there are many parameters, most of them are intended for solver development and have been pre-tuned. Here, we highlight only a few important parameters that are commonly adjusted in numerical experiments.

\begin{itemize}
\item \texttt{par.maxiter}: the maximum number of iterations allowed, default is $50000.$
\item \texttt{par.maxtime}: the maximum running time allowed, default is $3600$s.
\item \texttt{par.tol}: the tolerance on the KKT residual to terminate the algorithm, default is $10^{-6}.$
\item \texttt{par.verbose}: a flag indicating whether to print the solver's progress information during the run. It can be set to $0$ (silent mode) or $1$ (verbose mode, default).
\item \texttt{par.pfun}: the projection mapping for the constraint $X\in \P$ in (\ref{SDPG-multi}). It is a function handle with the input being the cell array of the variable: $X=\{X_1;\ldots;X_K\}$ and the output being the cell array of the variable after projection: 
$$
\{(\Pi_{\P}(X))_1;\ldots;(\Pi_{\P}(X))_K\}.
$$
\item \texttt{par.Ht} and \texttt{par.Hpfun}: the linear operator and projection mapping used to model the convex constraint, $\H(X)\in \Q$, in (\ref{SDPG-multi}). \texttt{Ht} has the same data structure as \texttt{At}. \texttt{Hpfun} is a function handle for the projection mapping onto the convex set $\Q.$ Both the input and output of \texttt{Hpfun} are $p$-dimensional vectors.
\item \texttt{par.Af}: a $K \times 1$ cell array containing the matrices $\{U_1; \ldots; U_K\}$ that define  ``facial" constraints of the form $X_i U_i = 0$ for the matrix blocks. In certain problems, such as the SDP relaxations of non-convex quadratic programs~\cite{burer2009copositive}, 
they
include constraints of the form $\langle X, UU^\top \rangle = 0$, which is equivalent to $XU = 0$. Motivated by the algorithms proposed in~\cite{hou2025low,hou2025rinnal+}, this kind of constraint is further equivalent to $U^\top R = 0$ in the low-rank setting. In SDPF+, when such constraints are present, the algorithm eliminates them by restricting the iterates to lie in the orthogonal complement of the column space of $U$.

\end{itemize} 

\end{enumerate}
The output \texttt{runhist} contains the following information:
\begin{itemize}
\item \texttt{runhist.fval}: the primal objective function value.
\item \texttt{runhist.dfval}: the dual objective function value. 
\item \texttt{runhist.ttime}: the total running time spent in SDPF+.
\item \texttt{runhist.iter}: the total iteration number of SDPF+. This is the iteration number of the feasible method plus the outer iteration number the of PALM.
\item \texttt{runhist.singular}: a flag indicating whether the problem is primal degenerate. It is set to $1$ if the problem is detected to be degenerate, and $0$ otherwise.
\item \texttt{runhist.numChol}: the number of Cholesky decompositions computed in the SDPF+.
\item \texttt{runhist.numCGiter}: the number of PCG iterations plus the number of ALM subiterations in SDPF+.
\item \texttt{runhist.pfeas, runhist.dfeas,runhist.comp}: the KKT residuals defined in (\ref{pfeas}) and (\ref{dfeascomp}).
\item \texttt{runhist.dfval}, \texttt{runhist.pdgap}: the dual objective function value and the duality gap, as defined in~(\ref{pdgap}). These values are available only when $f(\cdot)$ is a linear function and the fields \texttt{pfun}, \texttt{Ht}, and \texttt{Hpfun} are empty. Otherwise, they are set to \texttt{NaN}.
\item \texttt{runhist.X}: the primal solution. It is a \texttt{struct}-type array with a field \texttt{R}, which stores the matrix blocks in their low-rank form and the vector blocks in their original form.
\item \texttt{runhist.lam}: the Lagrange multipliers of the constraints $\A(X)=b,$ $l\leq\B(X)\leq u,$ $X\in \P$ and $\H(X)\in\Q$ respectively.
\item \texttt{runhist.lamk}: the Lagrange multiplier associated with the special constraint $X_{k,k} = I_k$. SDPF+ applies a specialized low-rank factorization to eliminate this constraint once it is detected in the problem. This technique has been employed in SDPF \cite{tang2024feasible}.
\item \texttt{runhist.S}: the dual slack matrix defined in (\ref{dSlack}). 
\end{itemize}

\section{Tables for numerical results}\label{App-table}

In this section, we list tables of some numerical experiment results.

\begin{center}
\begin{tiny}
\begin{longtable}{|c|ccccc|c|}
\caption{Test results for SDP problems of second-order reduced density matrices.}
\\
\hline
problem & algorithm & pfeas& dfeas & comp & fval & time \\ \hline
\endhead 
\label{test-2RDM}
AlH.1Sigma+.STO6G  & SDPF+     & 6.82e-06 & 8.71e-06 & 2.88e-06 & 2.4568419e+02 & 1.32e+03 \\
m=7230 & {\bf SDPNAL+}  & 8.66e-06 & 9.35e-06 & 1.65e-05 & 2.4568369e+02 & 6.78e+02 \\
\hline
B.2P.DZ  & SDPF+     & 9.90e-06 & 9.52e-06 & 1.81e-06 & 2.4573770e+01 & 1.46e+03 \\
m=7230 & {\bf SDPNAL+}  & 9.30e-06 & 5.56e-06 & 2.96e-06 & 2.4573585e+01 & 9.96e+02 \\
\hline
B.2P.SV  & SDPF+     & 6.57e-06 & 8.23e-06 & 2.66e-07 & 2.4559680e+01 & 4.67e+02 \\
m=4743 & SDPNAL+  & 6.81e-06 & 2.26e-06 & 4.94e-06 & 2.4559466e+01 & 4.52e+02 \\
\hline
B2.3Sigmag-.STO6G  & {\bf SDPF+}     & 4.48e-06 & 8.84e-06 & 1.60e-06 & 5.7336756e+01 & 9.24e+02 \\
m=7230 & SDPNAL+  & 9.85e-06 & 3.59e-06 & 5.69e-06 & 5.7335793e+01 & 1.24e+03 \\
\hline
BF.1Sigma+.STO6G  & SDPF+     & 4.43e-06 & 8.62e-06 & 1.81e-06 & 1.4247491e+02 & 1.21e+03 \\
m=7230 & {\bf SDPNAL+}  & 9.80e-06 & 3.90e-06 & 6.02e-06 & 1.4247334e+02 & 9.75e+02 \\
\hline
BH+.2Sigma+.STO6G  & SDPF+     & 3.86e-06 & 7.98e-06 & 2.71e-06 & 2.6980053e+01 & 2.00e+01 \\
m=948 & SDPNAL+  & 1.40e-06 & 3.58e-06 & 2.60e-06 & 2.6979873e+01 & 1.99e+01 \\
\hline
BH.1Sigma+.DZ  & SDPF+     & 1.51e-05 & 1.47e-05 & 1.04e-06 & 2.7336753e+01 & 3.62e+03 \\
m=15018 & {\bf SDPNAL+}  & 5.76e-06 & 5.70e-06 & 1.33e-05 & 2.7335536e+01 & 3.03e+03 \\
\hline
BN.3Pi.STO6G  & {\bf SDPF+}     & 2.71e-06 & 9.80e-06 & 8.51e-07 & 9.3285740e+01 & 1.51e+03 \\
m=7230 & SDPNAL+  & 8.81e-06 & 8.93e-06 & 8.04e-06 & 9.3285999e+01 & 1.62e+03 \\
\hline
BO.2Sigma+.STO6G  & {\bf SDPF+}     & 6.78e-06 & 7.18e-06 & 4.08e-06 & 1.1683343e+02 & 1.33e+03 \\
m=7230 & SDPNAL+  & 9.84e-06 & 4.12e-06 & 8.37e-06 & 1.1683227e+02 & 1.37e+03 \\
\hline
Be.1S.STO6G  & SDPF+     & 9.42e-06 & 4.71e-06 & 7.55e-07 & 1.4556158e+01 & 9.84e+00 \\
m=465 & {\bf SDPNAL+}  & 1.35e-06 & 4.02e-06 & 1.38e-06 & 1.4556168e+01 & 6.01e+00 \\
\hline
Be.1S.SV  & SDPF+     & 8.40e-06 & 6.57e-06 & 1.43e-07 & 1.4615922e+01 & 5.19e+02 \\
m=4743 & {\bf SDPNAL+}  & 6.78e-06 & 6.67e-06 & 6.58e-06 & 1.4615872e+01 & 2.66e+02 \\
\hline
Be2.1Sigmag+.infty.STO6G  & SDPF+     & 4.77e-06 & 8.42e-06 & 9.96e-07 & 2.9112395e+01 & 7.02e+02 \\
m=7230 & {\bf SDPNAL+}  & 1.37e-06 & 4.71e-06 & 4.85e-06 & 2.9112189e+01 & 3.65e+02 \\
\hline
BeF.2Sigma+.STO6G  & {\bf SDPF+}     & 4.22e-06 & 8.63e-06 & 2.78e-06 & 1.2763975e+02 & 1.16e+03 \\
m=7230 & SDPNAL+  & 9.99e-06 & 2.91e-06 & 1.23e-05 & 1.2763809e+02 & 1.20e+03 \\
\hline
BeH+.1Sigma+.STO6G  & {\bf SDPF+}     & 5.71e-06 & 9.00e-06 & 5.15e-06 & 1.6456779e+01 & 2.57e+01 \\
m=948 & SDPNAL+  & 2.70e-06 & 3.35e-06 & 2.92e-06 & 1.6456520e+01 & 2.63e+01 \\
\hline
BeH.2Sigma+.STO6G  & SDPF+     & 3.67e-06 & 8.70e-06 & 6.96e-07 & 1.6693095e+01 & 2.49e+01 \\
m=948 & {\bf SDPNAL+}  & 1.34e-06 & 3.99e-06 & 2.60e-06 & 1.6692939e+01 & 2.45e+01 \\
\hline
BeO.1Sigma+.STO6G  & SDPF+     & 4.89e-06 & 9.64e-06 & 7.54e-06 & 1.0192671e+02 & 1.31e+03 \\
m=7230 & {\bf SDPNAL+}  & 8.36e-06 & 8.52e-06 & 1.64e-05 & 1.0192578e+02 & 1.08e+03 \\
\hline
C.1S.STO6G  & SDPF+     & 2.66e-06 & 9.80e-06 & 3.83e-06 & 3.7520370e+01 & 1.77e+01 \\
m=465 & {\bf SDPNAL+}  & 6.45e-06 & 6.82e-06 & 1.08e-06 & 3.7520575e+01 & 6.29e+00 \\
\hline
C.3P.DZ  & {\bf SDPF+}     & 4.97e-06 & 9.09e-06 & 1.07e-05 & 3.7738361e+01 & 1.51e+03 \\
m=7230 & SDPNAL+  & 2.84e-13 & 3.36e-06 & 3.45e-05 & 3.7739141e+01 & 1.74e+03 \\
\hline
C.3P.STO6G  & SDPF+     & 2.32e-06 & 4.84e-06 & 2.58e-06 & 3.7592869e+01 & 9.22e+00 \\
m=465 & {\bf SDPNAL+}  & 2.40e-07 & 8.19e-06 & 8.74e-06 & 3.7592867e+01 & 5.89e+00 \\
\hline
C.3PSZ0.DZ  & SDPF+     & 4.08e-06 & 9.40e-06 & 3.26e-06 & 3.7739716e+01 & 1.91e+03 \\
m=7230 & {\bf SDPNAL+}  & 9.91e-06 & 3.57e-06 & 8.54e-07 & 3.7739725e+01 & 1.24e+03 \\
\hline
C2-.2Sigmag+.STO6G  & SDPF+     & 7.31e-06 & 9.05e-06 & 1.03e-06 & 9.0338358e+01 & 9.62e+02 \\
m=7230 & {\bf SDPNAL+}  & 7.64e-06 & 9.06e-06 & 2.99e-06 & 9.0336715e+01 & 7.72e+02 \\
\hline
C2.1Sigmag+.STO6G  & SDPF+     & 6.39e-06 & 9.25e-06 & 1.77e-06 & 9.0771933e+01 & 1.40e+03 \\
m=7230 & {\bf SDPNAL+}  & 9.73e-06 & 5.69e-06 & 1.02e-06 & 9.0769770e+01 & 1.01e+03 \\
\hline
CF.2Pir.STO6G  & SDPF+     & 8.64e-06 & 9.98e-06 & 2.98e-06 & 1.5914929e+02 & 1.47e+03 \\
m=7230 & {\bf SDPNAL+}  & 1.17e-12 & 4.13e-06 & 1.27e-05 & 1.5914776e+02 & 1.39e+03 \\
\hline
CH.2Pir.DZ  & {\bf SDPF+}     & 2.18e-05 & 1.12e-05 & 1.43e-05 & 4.1176874e+01 & 3.63e+03 \\
m=15018 & SDPNAL+  & 6.45e-06 & 8.74e-06 & 9.13e-06 & 4.1175943e+01 & 4.06e+03 \\
\hline
CH3+.1E''.STO6G  & SDPF+     & 4.94e-06 & 9.31e-06 & 7.70e-07 & 4.8751823e+01 & 2.28e+02 \\
m=2964 & {\bf SDPNAL+}  & 7.72e-06 & 7.54e-06 & 4.34e-06 & 4.8751264e+01 & 1.06e+02 \\
\hline
CH4.1A1.STO6G  & SDPF+     & 4.90e-06 & 6.68e-06 & 1.18e-06 & 5.3663795e+01 & 6.18e+02 \\
m=4743 & {\bf SDPNAL+}  & 5.59e-06 & 5.47e-06 & 1.70e-06 & 5.3663329e+01 & 2.11e+02 \\
\hline
CN.2Sigma+.STO6G  & SDPF+     & 5.20e-06 & 8.96e-06 & 6.21e-06 & 1.1104020e+02 & 1.32e+03 \\
m=7230 & {\bf SDPNAL+}  & 9.86e-06 & 1.12e-05 & 1.09e-05 & 1.1103988e+02 & 7.89e+02 \\
\hline
CO+.2Sigma+.STO6G  & SDPF+     & 7.86e-06 & 8.43e-06 & 2.14e-06 & 1.3481693e+02 & 1.35e+03 \\
m=7230 & {\bf SDPNAL+}  & 9.97e-06 & 6.84e-06 & 8.38e-06 & 1.3481557e+02 & 1.08e+03 \\
\hline
CO.1Sigma+.STO6G  & SDPF+     & 7.79e-06 & 8.55e-06 & 4.52e-06 & 1.3495777e+02 & 1.28e+03 \\
m=7230 & {\bf SDPNAL+}  & 3.92e-06 & 9.81e-06 & 4.25e-06 & 1.3495644e+02 & 7.69e+02 \\
\hline
FH2+.1A1.STO6G  & SDPF+     & 2.65e-06 & 9.71e-06 & 5.23e-06 & 1.0999157e+02 & 1.28e+02 \\
m=1743 & {\bf SDPNAL+}  & 2.30e-06 & 2.26e-06 & 1.18e-06 & 1.0999055e+02 & 6.34e+01 \\
\hline
H3.2A1'.DZ  & SDPF+     & 3.77e-06 & 9.51e-06 & 8.17e-07 & 3.3648050e+00 & 4.42e+01 \\
m=948 & {\bf SDPNAL+}  & 6.62e-07 & 3.07e-06 & 3.38e-07 & 3.3646782e+00 & 2.23e+01 \\
\hline
HF.1Sigma+.DZ  & {\bf SDPF+}     & 3.71e-05 & 1.95e-05 & 8.92e-05 & 1.0536495e+02 & 3.60e+03 \\
m=15018 & SDPNAL+  & 2.20e-05 & 7.38e-06 & 4.62e-05 & 1.0535898e+02 & 3.88e+03 \\
\hline
HLi2.2A1.STO6G  & {\bf SDPF+}     & 6.87e-06 & 7.91e-06 & 1.68e-06 & 1.9145482e+01 & 1.84e+03 \\
m=10593 & SDPNAL+  & 9.25e-08 & 1.60e-06 & 7.57e-06 & 1.9145128e+01 & 3.56e+03 \\
\hline
HN2+.1Sigma+.STO6G  & SDPF+     & 4.00e-06 & 9.55e-06 & 6.81e-06 & 1.3798216e+02 & 2.98e+03 \\
m=10593 & {\bf SDPNAL+}  & 9.36e-06 & 1.05e-05 & 3.17e-06 & 1.3797971e+02 & 1.37e+03 \\
\hline
HNO.1A'.STO6G  & SDPF+     & 2.27e-06 & 8.73e-06 & 5.24e-06 & 1.5967301e+02 & 2.86e+03 \\
m=10593 & {\bf SDPNAL+}  & 5.30e-06 & 1.20e-05 & 6.78e-06 & 1.5967075e+02 & 1.64e+03 \\
\hline
Li.2S.STO6G  & SDPF+     & 6.85e-06 & 8.30e-06 & 1.83e-06 & 7.4003136e+00 & 1.59e+01 \\
m=465 & {\bf SDPNAL+}  & 1.63e-06 & 9.82e-07 & 1.03e-06 & 7.4002434e+00 & 9.51e+00 \\
\hline
Li2.1Sigmag+.STO6G  & SDPF+     & 9.00e-06 & 4.15e-06 & 1.60e-06 & 1.6619854e+01 & 1.30e+03 \\
m=7230 & {\bf SDPNAL+}  & 5.04e-06 & 5.17e-06 & 6.71e-06 & 1.6619566e+01 & 5.82e+02 \\
\hline
LiF.1Sigma+.STO6G  & SDPF+     & 7.11e-06 & 7.45e-06 & 7.55e-07 & 1.1558150e+02 & 1.07e+03 \\
m=7230 & {\bf SDPNAL+}  & 9.86e-06 & 7.86e-06 & 5.58e-06 & 1.1558042e+02 & 8.99e+02 \\
\hline
LiH.1Sigma+.DZ  & SDPF+     & 5.50e-06 & 7.18e-06 & 4.11e-07 & 9.0038110e+00 & 3.56e+03 \\
m=15018 & {\bf SDPNAL+}  & 3.39e-06 & 6.01e-06 & 3.45e-06 & 9.0037725e+00 & 2.31e+03 \\
\hline
LiH.1Sigma+.STO6G  & SDPF+     & 8.57e-06 & 8.40e-06 & 4.31e-06 & 8.9673644e+00 & 3.75e+01 \\
m=948 & {\bf SDPNAL+}  & 4.73e-06 & 2.67e-06 & 1.71e-06 & 8.9672695e+00 & 2.95e+01 \\
\hline
LiOH.1Sigma+.STO6G  & SDPF+     & 3.38e-06 & 8.44e-06 & 6.02e-06 & 9.5671183e+01 & 2.29e+03 \\
m=10593 & {\bf SDPNAL+}  & 5.68e-06 & 6.45e-06 & 1.18e-05 & 9.5669720e+01 & 1.81e+03 \\
\hline
N.4S.DZ  & SDPF+     & 8.35e-06 & 4.18e-06 & 1.19e-06 & 5.4443861e+01 & 2.70e+03 \\
m=7230 & {\bf SDPNAL+}  & 7.37e-06 & 4.37e-06 & 1.64e-06 & 5.4444622e+01 & 7.88e+02 \\
\hline
N2+.2Sigmag+.STO6G  & SDPF+     & 6.70e-06 & 8.27e-06 & 8.60e-07 & 1.3145110e+02 & 1.23e+03 \\
m=7230 & {\bf SDPNAL+}  & 7.79e-06 & 6.41e-06 & 2.52e-06 & 1.3144894e+02 & 8.34e+02 \\
\hline
N2.1Sigmag+.STO6G  & SDPF+     & 6.93e-06 & 8.63e-06 & 6.62e-06 & 1.3232650e+02 & 1.39e+03 \\
m=7230 & {\bf SDPNAL+}  & 6.72e-07 & 1.11e-05 & 4.02e-07 & 1.3232514e+02 & 5.71e+02 \\
\hline
NH.1Delta.DZ  & {\bf SDPF+}     & 7.39e-06 & 1.91e-05 & 9.85e-06 & 5.8555132e+01 & 3.60e+03 \\
m=15018 & SDPNAL+  & 1.42e-05 & 7.61e-06 & 3.50e-05 & 5.8548842e+01 & 3.84e+03 \\
\hline
NH.3Sigma-.DZ  & {\bf SDPF+}     & 4.45e-05 & 1.75e-04 & 4.35e-04 & 5.8617473e+01 & 3.65e+03 \\
m=15018 & SDPNAL+  & 1.43e-05 & 5.96e-06 & 1.83e-05 & 5.8614018e+01 & 3.93e+03 \\
\hline
NH2-.1A1.STO6G  & SDPF+     & 6.56e-06 & 9.67e-06 & 1.16e-06 & 6.2707265e+01 & 1.11e+02 \\
m=1743 & {\bf SDPNAL+}  & 3.69e-06 & 3.77e-06 & 1.47e-06 & 6.2706412e+01 & 5.81e+01 \\
\hline
NH3+.2A2".STO6G  & SDPF+     & 4.57e-06 & 9.99e-06 & 1.62e-06 & 6.7506589e+01 & 2.49e+02 \\
m=2964 & {\bf SDPNAL+}  & 5.70e-06 & 5.03e-06 & 3.86e-06 & 6.7506042e+01 & 1.18e+02 \\
\hline
NH4+.1A1.STO6G  & SDPF+     & 8.86e-06 & 9.59e-06 & 1.37e-07 & 7.2777720e+01 & 6.19e+02 \\
m=4743 & {\bf SDPNAL+}  & 3.60e-06 & 9.64e-06 & 1.62e-05 & 7.2778616e+01 & 3.04e+02 \\
\hline
Na.2S.STO6G  & SDPF+     & 9.33e-06 & 8.13e-06 & 2.30e-06 & 1.6107858e+02 & 5.34e+02 \\
m=4743 & {\bf SDPNAL+}  & 3.48e-06 & 6.17e-06 & 3.62e-05 & 1.6107816e+02 & 3.11e+02 \\
\hline
NaH.1Sigma+.STO6G  & SDPF+     & 4.64e-06 & 9.81e-06 & 2.17e-06 & 1.6482460e+02 & 1.18e+03 \\
m=7230 & {\bf SDPNAL+}  & 9.80e-06 & 9.42e-06 & 2.49e-05 & 1.6482336e+02 & 6.06e+02 \\
\hline
Ne.1S.DZ  & SDPF+     & 6.07e-06 & 8.94e-06 & 7.86e-06 & 1.2864466e+02 & 1.72e+03 \\
m=7230 & {\bf SDPNAL+}  & 9.15e-07 & 5.10e-06 & 5.45e-05 & 1.2864187e+02 & 5.86e+02 \\
\hline
O.1D.DZ  & SDPF+     & 9.53e-06 & 9.19e-06 & 4.15e-06 & 7.4791985e+01 & 1.52e+03 \\
m=7230 & {\bf SDPNAL+}  & 7.11e-06 & 9.20e-06 & 3.04e-05 & 7.4791127e+01 & 6.53e+02 \\
\hline
O.3P.DZ  & SDPF+     & 8.03e-06 & 9.56e-06 & 3.24e-06 & 7.4873477e+01 & 1.69e+03 \\
m=7230 & {\bf SDPNAL+}  & 6.27e-07 & 7.84e-06 & 1.99e-05 & 7.4872796e+01 & 7.02e+02 \\
\hline
O.3PSZ0.DZ  & SDPF+     & 3.13e-06 & 5.00e-06 & 5.08e-06 & 7.4873116e+01 & 1.91e+03 \\
m=7230 & {\bf SDPNAL+}  & 7.50e-06 & 7.48e-06 & 2.01e-06 & 7.4874502e+01 & 8.63e+02 \\
\hline
O2+.2Pig.STO6G& SDPF+ & 7.29e-06 & 9.93e-06 & 2.59e-06 & 1.7913508e+02& 1.96e+03 \\
m=7230 & {\bf SDPNAL+} & 5.86e-06 & 8.66e-06 & 2.94e-06 & 1.7913201e+02& 6.84e+02 \\
\hline
\end{longtable}
\end{tiny}
\end{center}

\begin{center}
\begin{tiny}
\begin{longtable}{|c|ccccc|c|}
\caption{Test results for SDP problems from Hans Mittelmann benchmarks. 
Here comp' denotes $\min\{{\rm comp,pdgap}\}$.}
\\
\hline
problem & algorithm & pfeas& dfeas & comp' & fval & time \\ \hline
\endhead
\label{test-HMT}
Alh\_1.r20 &  SDPF+ & 1.45e-06 & 4.35e-06 & 2.04e-06 & 2.4568341e+02& 3.68e+03 \\
nblk=22& {\bf SDPNAL+} & 2.82e-07 & 9.69e-07 & 7.63e-06 & 2.4568507e+02& 2.06e+03 \\
m=7230& SDPT3 & -& - & - & -& - \\
\hline
BH2.r14 & SDPF+ & 6.78e-07 & 9.74e-07 & 5.98e-08 & 3.0430316e+01& 7.60e+02 \\
nblk=22 & {\bf SDPNAL+} & 2.31e-07 & 1.17e-06 & 2.52e-06 & 3.0430377e+01& 1.87e+02 \\
m=1743 & SDPT3 & 1.71e-09 & 0.00e+00 & 5.31e-07 & 3.0430147e+01& 3.99e+02 \\
\hline
Bex2\_1\_5 & SDPF+ & 5.10e-07 & 9.66e-07 & 2.94e-09 & 2.9779301e+00& 2.04e+02 \\
nblk=33 & {\bf SDPNAL+} & 5.18e-07 & 1.02e-08 & 1.13e-10 & 2.9779225e+00& 3.25e+01 \\
m=3002 & SDPT3 & 6.22e-05 & 0.00e+00 & 1.80e-06 & 2.9778612e+00& 4.16e+02 \\
\hline
Bst\_jcbpaf2 & SDPF+ & 3.61e-07 & 6.23e-07 & 2.87e-09 & 7.9487455e-02& 2.55e+02 \\
nblk=35 & {\bf SDPNAL+} & 3.19e-07 & 3.46e-08 & 1.51e-08 & 7.9483190e-02& 4.17e+01 \\
m=3002 & SDPT3 & 1.30e-07 & 0.00e+00 & 1.18e-08 & 7.9485765e-02& 4.66e+02 \\
\hline
CH2\_1.r14 & SDPF+ & 7.01e-07 & 9.10e-07 & 2.30e-07 & 4.4853922e+01& 6.35e+02 \\
nblk=22 & {\bf SDPNAL+} & 9.03e-07 & 1.23e-06 & 2.32e-06 & 4.4853953e+01& 1.52e+02 \\
m=1743 & SDPT3 & 8.45e-10 & 0.00e+00 & 4.59e-07 & 4.4853801e+01& 3.87e+02 \\
\hline
G40\_mb & {\bf SDPF+} & 1.00e-08 & 2.63e-08 & 7.28e-10 & -2.8643232e+03& 1.18e+01 \\
nblk= 1 & SDPNAL+ & 3.99e-07 & 7.52e-11 & 5.55e-07 & -2.8643297e+03& 3.60e+03 \\
m=2001 & SDPT3 & 1.96e-12 & 0.00e+00 & 6.28e-07 & -2.8643196e+03& 3.80e+01 \\
\hline
G40mc & {\bf SDPF+} & 8.14e-09 & 1.95e-08 & 2.90e-10 & -5.7295791e+03& 7.20e+00 \\
nblk= 1 & SDPNAL+ & 3.19e-07 & 8.35e-08 & 3.86e-08 & -5.7295795e+03& 2.05e+02 \\
m=2000 & SDPT3 & 1.40e-11 & 0.00e+00 & 1.93e-07 & -5.7295769e+03& 3.45e+01 \\
\hline
G48\_mb & {\bf SDPF+} & 2.22e-09 & 0.00e+00 & 5.52e-07 & 8.2171572e+00& 1.07e+01 \\
nblk= 1 & SDPNAL+ & 5.56e-07 & 1.26e-09 & 2.94e-05 & 8.2191558e+00& 1.02e+03 \\
m=3001 & SDPT3 & 2.31e-12 & 0.00e+00 & 1.26e-07 & 8.2171648e+00& 7.90e+01 \\
\hline
G48mc & {\bf SDPF+} & 7.53e-12 & 0.00e+00 & 4.82e-08 & -1.1999999e+04& 6.85e+00 \\
nblk= 1 & SDPNAL+ & 1.55e-07 & 7.85e-11 & 1.54e-08 & -1.2000000e+04& 9.17e+02 \\
m=3000 & SDPT3 & 1.22e-15 & 0.00e+00 & 1.73e-08 & -1.2000000e+04& 4.19e+01 \\
\hline
G55mc & {\bf SDPF+} & 1.06e-09 & 2.27e-08 & 6.44e-08 & -2.2078920e+04& 2.26e+01 \\
nblk= 1 & SDPNAL+ & 4.55e-10 & 1.09e-05 & 6.25e-06 & -2.2078888e+04& 3.61e+03 \\
 m=5000 & SDPT3 & 3.45e-14 & 0.00e+00 & 2.33e-07 & -2.2078911e+04& 3.09e+02 \\
\hline
G59mc & {\bf SDPF+} & 6.55e-09 & 4.61e-08 & 1.36e-08 & -1.4624653e+04& 3.01e+01 \\
nblk= 1 & SDPNAL+ & 6.63e-07 & 1.53e-06 & 5.46e-06 & -1.4624654e+04& 3.66e+03 \\
m=5000 & SDPT3 & 4.13e-12 & 0.00e+00 & 2.04e-07 & -1.4624647e+04& 4.24e+02 \\
\hline
G60\_mb & {\bf SDPF+} & 4.14e-09 & 1.09e-08 & 1.29e-07 & 1.9280682e+03& 5.64e+01 \\
nblk= 1 & SDPNAL+ & - & - & - & -& - \\
m=7001 & SDPT3 & 3.05e-13 & 1.65e-18 & 2.53e-07 & 1.9280690e+03& 7.71e+02 \\
\hline
 G60mc & {\bf SDPF+} & 3.16e-09 & 1.00e-08 & 4.78e-08 & -3.0444536e+04& 5.49e+01 \\
nblk= 1 & SDPNAL+ & - & - & - & -& - \\
m=7000 & SDPT3 & 1.89e-14 & 0.00e+00 & 5.35e-07 & -3.0444504e+04& 7.85e+02 \\
\hline
H3O+\_.r16 & SDPF+ & 9.11e-07 & 9.83e-07 & 7.93e-07 & 9.0107034e+01& 1.97e+03 \\
nblk=22 & {\bf SDPNAL+} & 1.27e-07 & 1.38e-06 & 1.62e-06 & 9.0106969e+01& 2.55e+02 \\
m=2964 & SDPT3 & 1.03e-11 & 0.00e+00 & 6.70e-07 & 9.0106681e+01& 1.87e+03 \\
\hline
NH2-.r14 & SDPF+ & 3.49e-07 & 9.99e-07 & 4.10e-07 & 6.2706409e+01& 1.02e+03 \\
nblk=22 & {\bf SDPNAL+} & 4.14e-07 & 1.23e-06 & 8.43e-07 & 6.2706337e+01& 8.62e+01 \\
m=1743 & SDPT3 & 3.11e-08 & 4.98e-07 & 6.24e-07 & 6.2706343e+01& 5.46e+02 \\
\hline
NH3\_.r16 & SDPF+ & 6.41e-07 & 9.62e-07 & 6.87e-07 & 6.7925265e+01& 2.62e+03 \\
 nblk=22 & {\bf SDPNAL+} & 3.42e-07 & 1.39e-06 & 1.91e-06 & 6.7925297e+01& 2.49e+02 \\
m=2964 & SDPT3 & 1.13e-11 & 0.00e+00 & 4.00e-07 & 6.7924925e+01& 1.87e+03 \\
\hline
NH4+.r18 & SDPF+ & 2.81e-05 & 5.60e-05 & 1.43e-05 & 7.2777756e+01& 3.63e+03 \\
nblk=22 & {\bf SDPNAL+} & 9.72e-07 & 9.42e-07 & 2.54e-06 & 7.2776721e+01& 9.33e+02 \\
m=4743 & SDPT3 & - & - & - & -& - \\
\hline
biggs & SDPF+ & - & - & - &    -& - \\
nblk=48 & SDPNAL+ & - & - & - & -& - \\
m=1819 & {\bf SDPT3} & 1.31e-09 & 0.00e+00 & 2.90e-07 & -1.4142583e+03& 3.21e+02 \\
\hline
 broyden25 & SDPF+ & 3.14e-07 & 8.03e-10 & 2.18e-07 & 3.5536241e+01& 1.46e+03 \\
nblk=19 & SDPNAL+ & - & - & - & -& - \\
m=15971 & {\bf SDPT3} & 2.10e-10 & 0.00e+00 & 8.41e-07 & 3.5536255e+01& 1.10e+03 \\
\hline
buck3 & SDPF+ & 4.94e-07 & 1.42e-12 & 7.10e-07 & -6.0760518e+02& 9.35e+02 \\
nblk= 3 & SDPNAL+ & - & - & - & -& - \\
m=544 & {\bf SDPT3} & 2.41e-06 & 0.00e+00 & 3.13e-06 & -6.0760227e+02& 5.36e+00 \\
\hline
 buck4 & SDPF+ & 9.89e-07 & 3.75e-08 & 7.83e-08 & -4.8613148e+02& 5.78e+02 \\
 nblk= 3 & SDPNAL+ & - & - & - & - & - \\
m=1200 & {\bf SDPT3} & 2.57e-07 & 0.00e+00 & 7.34e-07 & -4.8614176e+02& 1.96e+01 \\
\hline
buck5 & SDPF+ & 9.95e-07 & 6.10e-08 & 8.37e-07 & -4.3616210e+02& 2.22e+03 \\
nblk= 3 & SDPNAL+ & - & - & - & -& - \\
m=3280 & {\bf SDPT3} & 1.49e-06 & 0.00e+00 & 7.63e-05 & -4.3619638e+02& 1.85e+02 \\
\hline
cancer& SDPF+ & 2.39e-06 & 3.17e-06 & 6.18e-08 & -2.7623801e+04 & 3.60e+03 \\
nblk= 1 & SDPNAL+ & 9.93e-15 & 2.82e-11 & 1.89e-05 & -2.7625376e+04& 3.60e+03 \\
m=10469 & {\bf SDPT3} & 9.26e-08 & 0.00e+00 & 2.89e-04 & -2.7607459e+04& 8.59e+01 \\
\hline
checker& {\bf SDPF+} & 4.66e-09 & 9.69e-11 & 2.87e-07 & 3.3038867e+03& 4.54e+01 \\
nblk= 1 & SDPNAL+ & 1.28e-14 & 2.21e-06 & 3.33e-04 & 3.3065285e+03& 3.61e+03 \\
m=3970 & SDPT3 & 3.13e-11 & 0.00e+00 & 2.28e-07 & 3.3038861e+03& 1.23e+02 \\
\hline
cnhil10 & SDPF+ & 8.13e-07 & 9.37e-07 & 2.51e-10 & 0.0000000e+00& 1.68e+02 \\
nblk= 1 & SDPNAL+ & 8.94e-08 & 2.23e-07 & 3.04e-06 & 0.0000000e+00& 1.95e+01 \\
m=5005 & {\bf SDPT3} & 3.45e-08 & 0.00e+00 & 2.73e-04 & 0.0000000e+00& 1.70e+01 \\
\hline
cnhil8 & SDPF+ & 9.46e-07 & 7.37e-07 & 1.60e-09 & 0.0000000e+00& 1.72e+01 \\
nblk= 1 & SDPNAL+ & 2.70e-07 & 3.13e-07 & 7.95e-07 & 0.0000000e+00& 9.04e+00 \\
m=1716 & {\bf SDPT3} & 1.72e-08 & 0.00e+00 & 6.77e-07 & 0.0000000e+00& 1.90e+00 \\
\hline
 cphil10& SDPF+ & 1.73e-13 & 1.29e-10 & 3.64e-13 & 0.0000000e+00& 9.00e+00 \\
nblk= 1 & {\bf SDPNAL+} & 1.59e-16 & 9.13e-17 & 1.52e-16 & 0.0000000e+00& 2.80e-01 \\
m=5005 & SDPT3 & 1.93e-14 & 0.00e+00 & 2.04e-08 & 0.0000000e+00& 9.87e+00 \\
\hline
cphil12 & SDPF+ & 3.09e-13 & 4.69e-10 & 6.70e-13 & 0.0000000e+00& 5.70e+01 \\
nblk= 1 & {\bf SDPNAL+} & 1.63e-16 & 1.26e-16 & 1.59e-16 & 0.0000000e+00& 4.16e-01 \\
m=12376 & SDPT3 & 1.48e-13 & 0.00e+00 & 2.01e-07 & 0.0000000e+00& 6.47e+01 \\
\hline
dia\_patch & {\bf SDPF+} & 2.28e-07 & 3.30e-07 & 4.83e-07 & 5.4174249e-06& 1.53e+03 \\
nblk= 1 & SDPNAL+ & - & - & - & -& - \\
m=5478 & SDPT3 & - & - & - & -& - \\
\hline
e*quad* & SDPF+ & 2.93e-05 & 1.96e-04 & 1.43e-04 & 4.9544525e+03& 3.60e+03 \\
nblk= 3 & SDPNAL+ & 1.20e-09 & 1.72e-06 & 1.73e-04 & 4.9270250e+03& 3.22e+02 \\
m=5984 & {\bf SDPT3} & 1.43e-07 & 3.87e-08 & 8.90e-07 & 4.9273596e+03& 9.21e+01 \\
\hline
e*stable* & SDPF+ & 5.59e-07 & 1.35e-07 & 9.69e-07 & -1.9843376e-01& 1.29e+03 \\
nblk=19 & SDPNAL+ & 9.20e-07 & 4.58e-07 & 1.38e-05 & -1.9846413e-01& 1.16e+02 \\
m=5984 & {\bf SDPT3} & 1.05e-09 & 6.08e-08 & 9.34e-06 & -1.9843397e-01& 8.21e+01 \\
\hline
fap09 & {\bf SDPF+} & 9.99e-07 & 8.68e-07 & 9.57e-12 & 1.0780819e+01& 2.81e+02 \\
nblk= 2 & SDPNAL+ & 2.74e-07 & 4.27e-04 & 2.15e-06 & 1.0809393e+01& 4.74e+00 \\
m=30276 & SDPT3 & - & - & - & -& - \\
\hline
 foot & SDPF+ & 5.74e-09 & 0.00e+00 & 1.36e-08 & -5.8529815e+05& 6.14e+01 \\
nblk= 1 & SDPNAL+ & 3.09e-07 & 1.89e-10 & 1.79e-05 & -5.8534088e+05& 3.60e+03 \\
m=2209 & {\bf SDPT3} & 3.97e-06 & 0.00e+00 & 8.20e-06 & -5.8528857e+05& 4.73e+01 \\
\hline
hand & {\bf SDPF+} & 1.84e-08 & 2.56e-08 & 3.90e-09 & -2.4747779e+04& 6.84e+00 \\
nblk= 1 & SDPNAL+ & 7.72e-07 & 6.20e-10 & 1.24e-07 & -2.4747790e+04& 1.82e+03 \\
m=1297 & SDPT3 & 1.55e-08 & 0.00e+00 & 5.53e-07 & -2.4747752e+04& 1.26e+01 \\
\hline
ice\_2.0 & {\bf SDPF+} & 2.43e-09 & 1.63e-10 & 4.08e-07 & 6.8083911e+03& 2.14e+02 \\
nblk= 1 & SDPNAL+ & - & - & - & - & - \\
m=8113 & SDPT3 & 3.87e-12 & 0.00e+00 & 3.80e-07 & 6.8083914e+03& 7.79e+02 \\
\hline
inc\_1200 & SDPF+ & 2.32e-06 & 4.75e-07 & 8.83e-10 & -1.1575928e+00& 3.60e+03 \\
nblk= 2 & SDPNAL+ & 3.79e-12 & 8.09e-06 & 1.90e-03 & 8.3700086e-02& 3.60e+03 \\
m=5175 & {\bf SDPT3} & 4.81e-08 & 0.00e+00 & 1.23e-05 & -1.1566568e+00& 6.57e+01 \\
\hline
inc\_600 & SDPF+ & 9.84e-07 & 2.56e-07 & 4.09e-09 & -6.6855340e-01& 2.46e+03 \\
nblk= 2 & SDPNAL+ & 3.65e-14 & 6.75e-06 & 1.81e-03 & -2.3678686e-02& 3.60e+03 \\
m=2515 & {\bf SDPT3} & 2.74e-08 & 2.57e-18 & 1.94e-06 & -6.6843406e-01& 1.60e+01 \\
\hline
neosfbr25 & SDPF+ & 2.66e-05 & 1.13e-05 & 4.74e-07 & 3.4537603e+03& 3.60e+03 \\
nblk= 1 & SDPNAL+ & 1.41e-09 & 1.81e-07 & 2.78e-06 & 3.4535482e+03& 2.50e+03 \\
m=14376 & {\bf SDPT3} & 2.58e-11 & 0.00e+00 & 8.31e-08 & 3.4536139e+03& 7.76e+02 \\
\hline
neosfbr30e8 & SDPF+ & 1.26e-04 & 1.05e-04 & 4.51e-04 & 5.7047831e+03& 3.63e+03 \\
nblk= 1 & SDPNAL+ & 3.00e-10 & 1.29e-06 & 1.85e-05 & 5.6961603e+03& 3.60e+03 \\
m=25201 & {\bf SDPT3} & 1.98e-12 & 0.00e+00 & 7.35e-08 & 5.6954664e+03& 2.82e+03 \\
\hline
neu1 & SDPF+ & 3.46e-07 & 6.62e-07 & 7.37e-08 & -5.0331671e-06& 4.13e+02 \\
nblk= 2 & SDPNAL+ & 1.05e-13 & 2.11e-09 & 1.88e-04 & -4.6682975e-04& 5.09e+02 \\
m=3003 & {\bf SDPT3} & 8.32e-09 & 9.23e-09 & 9.12e-05 & -1.0776363e-05& 1.27e+02 \\
\hline
neu1g & SDPF+ & 6.20e-07 & 4.45e-07 & 9.61e-09 & 1.2500000e+02& 1.12e+02 \\
nblk= 1 & SDPNAL+ & 1.62e-13 & 1.79e-09 & 1.80e-06 & 1.2499954e+02& 4.70e+02 \\
m=3002 & {\bf SDPT3} & 4.19e-10 & 0.00e+00 & 1.99e-07 & 1.2500000e+02& 1.04e+02 \\
\hline
neu2 & SDPF+ & - & - & - & -& - \\
nblk= 2 & SDPNAL+ & - & - & - & -& - \\
m=3003 & {\bf SDPT3} & 5.05e-08 & 1.24e-07 & 2.12e-04 & -6.0403737e+02& 1.26e+02 \\
\hline
neu2c & SDPF+ & - & - & - & -& - \\
nblk=13 & SDPNAL+ & - & - & - & -& - \\
m=3002 & {\bf SDPT3} & 3.98e-05 & 0.00e+00 & 9.13e-03 & 3.4472963e+04& 3.39e+02 \\
\hline
neu2g & SDPF+ & 1.25e-03 & 2.85e-08 & 1.46e-03 & 3.3819730e+04& 3.60e+03 \\
nblk= 1 & SDPNAL+ & - & - & - & -& - \\
m=3002 & {\bf SDPT3} & 3.06e-08 & 0.00e+00 & 6.80e-07 & 3.4100002e+04& 1.02e+02 \\
\hline
neu3g & SDPF+ & 6.94e-07 & 2.89e-07 & 9.76e-07 & 3.1254591e-07& 6.88e+02 \\
 nblk= 1 & {\bf SDPNAL+} & 9.28e-07 & 7.95e-08 & 5.03e-07 & 6.4084985e-07& 6.52e+01 \\
m=8007 & SDPT3 & 6.53e-15 & 0.00e+00 & 8.32e-05 & 1.3150173e-05& 1.96e+02 \\
\hline
p\_auss2\_3.0 & {\bf SDPF+} & 2.44e-09 & 8.92e-11 & 8.07e-07 & 8.6181059e+03& 2.67e+02 \\
nblk=1 & SDPNAL+ & - & - & - & -& - \\
m=9115 & SDPT3 & 1.44e-12 & 0.00e+00 & 8.54e-07 & 8.6181093e+03& 9.47e+02 \\
\hline
rabmo& SDPF+ & 1.84e-05 & 7.47e-06 & 3.27e-06 & -3.7182372e+00& 3.60e+03 \\
nblk= 2 & SDPNAL+ & 6.04e-13 & 5.61e-07 & 2.44e-07 & -3.7274501e+00& 1.60e+02 \\
m=5004 & {\bf SDPT3} & 1.56e-06 & 3.66e-06 & 2.71e-04 & -3.7272538e+00& 1.04e+02 \\
\hline
reimer5  & SDPF+ & 3.77e-03 & 1.64e-04 & 4.33e-04 & -1.5163447e+01& 3.60e+03 \\
nblk= 2 & SDPNAL+ & - & - & - & -& - \\
m=6187 & {\bf SDPT3} & 8.11e-06 & 1.03e-07 & 1.40e-05 & -1.5183781e+01& 2.71e+03 \\
\hline
 r1\_2000 & SDPF+ & 3.49e-09 & 1.40e-07 & 2.34e-05 & -4.8997334e+05& 3.60e+03 \\
nblk= 1 & SDPNAL+ & 1.21e-04 & 9.09e-12 & 2.71e-07 & -4.9029307e+05& 3.60e+03 \\
m=2001 & {\bf SDPT3} & 1.04e-08 & 0.00e+00 & 3.93e-07 & -4.9029215e+05& 3.52e+01 \\
\hline
 rose13 & SDPF+ & 8.89e-07 & 6.16e-07 & 5.60e-07 & 1.1999921e+01& 8.05e+01 \\
 nblk= 1 & SDPNAL+ & 1.84e-15 & 8.86e-09 & 6.13e-05 & 1.1998465e+01& 1.97e+02 \\
m=2379 & {\bf SDPT3} & 8.38e-11 & 0.00e+00 & 3.72e-07 & 1.2000000e+01& 7.83e+00 \\
\hline
rose15 & SDPF+ & 8.40e-07 & 1.52e-07 & 9.37e-07 & -8.1462754e-05& 2.84e+03 \\
nblk= 2 & SDPNAL+ & 4.59e-16 & 6.76e-09 & 8.53e-04 & -1.9993087e-03& 3.56e+02 \\
m=3860 & {\bf SDPT3} & 3.14e-09 & 4.16e-08 & 1.44e-03 & -6.3462959e-06&3.46e+01 \\
\hline
sensor\_1000b  & SDPF+ & 9.68e-07 & 8.62e-07 & 2.60e-09 & 2.3355746e+00& 2.42e+02 \\
nblk= 2 & SDPNAL+ & 4.08e-07 & 1.35e-07 & 1.87e-05 & 2.3347799e+00& 1.38e+03 \\
m=5549 & {\bf SDPT3} & 5.60e-10 & 0.00e+00 & 6.15e-08 & 2.3355818e+00& 6.28e+01 \\
\hline
sensor\_500b  & SDPF+ & 8.63e-07 & 1.57e-07 & 2.32e-08 & 1.7389448e+01& 6.67e+01 \\
nblk= 2 & SDPNAL+ & 9.10e-07 & 1.32e-08 & 2.30e-09 & 1.7389458e+01& 9.99e+01 \\
m=3540 & {\bf SDPT3} & 7.40e-11 & 0.00e+00 & 1.12e-07 & 1.7389454e+01& 2.26e+01 \\
\hline
shmup3& SDPF+ & 4.07e-07 & 6.66e-08 & 9.88e-07 & -2.0988377e+03& 5.27e+02 \\
nblk= 3 & SDPNAL+ & - & - & - & -& - \\
m=420 & {\bf SDPT3} & 2.06e-09 & 0.00e+00 & 7.23e-07 & -2.0988351e+03& 3.53e+01 \\
\hline
shmup4 & SDPF+ & 6.37e-07 & 5.68e-08 & 9.88e-07 & -7.9925517e+03& 1.83e+03 \\
nblk= 3 & SDPNAL+ & - & - & - & -& - \\
m=800 & {\bf SDPT3} & 6.91e-09 & 1.89e-14 & 3.25e-07 & -7.9925465e+03& 1.80e+02 \\
\hline
shmup5  & SDPF+ & 2.30e-05 & 3.23e-07 & 8.69e-05 & -2.3859149e+04& 3.60e+03 \\
nblk= 3 & SDPNAL+ & - & - & - & -& - \\
m=1800 & {\bf SDPT3} & 4.58e-08 & 0.00e+00 & 5.11e-07 & -2.3858904e+04& 1.14e+03 \\
\hline
 spar060 & {\bf SDPF+} & 4.19e-07 & 0.00e+00 & 1.89e-09 & -1.2780000e+03& 4.65e+00 \\
nblk= 2 & SDPNAL+ & 6.41e-07 & 1.76e-12 & 6.14e-08 & -1.2779990e+03& 1.99e+01 \\
m=7261 & SDPT3 & 4.04e-12 & 0.00e+00 & 6.40e-08 & -1.2779999e+03& 7.36e+01 \\
\hline
swissroll & {\bf SDPF+} & 7.46e-04 & 8.92e-10 & 3.25e-03 & -5.5490868e+05& 3.60e+03 \\
 nblk= 1 & SDPNAL+ & - & - & - & -& - \\
m=3380 & SDPT3 & - & -& - & -& - \\
\hline
taha1a & SDPF+ & 1.38e-07 & 8.89e-07 & 9.70e-07 & -9.9999539e-01& 1.83e+02 \\
nblk=14 & {\bf SDPNAL+} & 4.89e-07 & 9.25e-08 & 2.65e-08 & -9.9999872e-01& 5.85e+01 \\
m=3002 & SDPT3 & 1.26e-05 & 8.84e-18 & 8.69e-07 & -1.0000070e+00& 2.32e+02 \\
\hline
 taha1b & SDPF+ & 5.00e-07 & 2.00e-07 & 1.66e-06 & -7.7329713e-01& 3.60e+03 \\
nblk=22 & SDPNAL+ & 1.00e-13 & 3.85e-08 & 4.08e-05 & -7.7352821e-01& 3.60e+03 \\
m=8007 & {\bf SDPT3} & 4.64e-12 & 0.00e+00 & 7.16e-06 & -7.7328713e-01& 2.77e+02 \\
\hline
taha1c & SDPF+ & 6.25e-07 & 8.71e-07 & 2.39e-07 & -9.9999795e-01& 7.87e+02 \\
nblk=14 & {\bf SDPNAL+} & 2.73e-12 & 1.79e-07 & 1.58e-06 & -1.0000046e+00& 3.41e+02 \\
m=6187 & SDPT3 & 9.51e-03 & 1.06e-17 & 2.02e-03 & -1.0000122e+00& 2.44e+03 \\
\hline
theta12 & SDPF+ & 5.89e-07 & 4.07e-08 & 6.92e-10 & -9.2801485e+01& 2.30e+01 \\
nblk= 1 & {\bf SDPNAL+} & 9.24e-07 & 7.80e-09 & 2.59e-10 & -9.2801694e+01& 1.74e+01 \\
m=17979 & SDPT3 & 1.28e-12 & 0.00e+00 & 2.72e-07 & -9.2801537e+01& 1.86e+02 \\
\hline
theta123 & SDPF+ & 3.72e-07 & 1.56e-07 & 1.60e-12 & -2.4668675e+01& 3.83e+02 \\
nblk= 1 & {\bf SDPNAL+} & 5.74e-07 & 2.76e-08 & 7.42e-11 & -2.4668640e+01& 1.61e+01 \\
m=90020 & SDPT3 & - & - & - & -& - \\
\hline
 t\_texture & SDPF+ & 1.31e-04 & 7.23e-10 & 6.85e-04 &  3.4445867e+02& 3.60e+03 \\
nblk= 1 & SDPNAL+ & 1.23e-07 & 3.52e-09 & 8.79e-04 & 3.4410041e+02& 3.60e+03 \\
m=1802 & {\bf SDPT3} & 4.67e-08 & 0.00e+00 & 3.53e-04 & 3.4460806e+02& 4.44e+01 \\
\hline
trto3 & SDPF+ & 8.21e-07 & 3.86e-12 & 5.01e-07 & -1.2799974e+04& 1.24e+02 \\
nblk= 2 & SDPNAL+ & - & - & - & -& - \\
m=544 & {\bf SDPT3} & 2.32e-05 & 0.00e+00 & 9.32e-07 & -1.2799982e+04& 2.11e+00 \\
\hline
trto4 & SDPF+ & 4.51e-07 & 6.26e-17 & 4.49e-07 & -1.2765821e+04& 1.11e+03 \\
nblk= 2 & SDPNAL+ & - & - & - & -& - \\ 
m=1200 & {\bf SDPT3} & 7.52e-06 & 0.00e+00 & 1.08e-06 & -1.2765805e+04& 9.69e+00 \\
\hline
 trto5 & SDPF+ & 5.01e-04 & 1.57e-10 & 2.13e-03 & -1.2808388e+04& 3.60e+03 \\
nblk= 2 & SDPNAL+ & - & - & - & -& - \\ 
m=3280 & {\bf SDPT3} & 8.39e-05 & 0.00e+00 & 1.41e-05 & -1.2799641e+04& 7.08e+01 \\
\hline
vibra3 & SDPF+ & 3.70e-07 & 6.34e-12 & 9.90e-07 & -1.7261308e+02& 6.98e+02 \\
nblk= 3 & SDPNAL+ & - & - & - & -& - \\ 
m=544 & {\bf SDPT3} & 3.11e-06 & 0.00e+00 & 1.19e-06 & -1.7261270e+02& 4.59e+00 \\
\hline
 vibra4 & SDPF+ & 1.16e-03 & 1.59e-09 & 5.54e-03 & -1.6547537e+02& 3.60e+03 \\
nblk= 3 & SDPNAL+ & - & - & - & -& - \\ 
m=1200 & {\bf SDPT3} & 9.61e-07 & 0.00e+00 & 9.72e-07 & -1.6561301e+02& 2.66e+01 \\
\hline
 vibra5 & SDPF+ & 5.10e-04 & 4.78e-07 & 1.49e-03 & -1.2819633e+02& 3.60e+03 \\
nblk= 3 & SDPNAL+ & - & - & - & -& - \\ 
m=3280 & {\bf SDPT3} & 9.14e-06 & 0.00e+00 & 2.58e-05 & -1.6589585e+02& 2.91e+02 \\
\hline
yalsdp & SDPF+ & 4.88e-07 & 6.88e-11 & 8.80e-08 & -1.7921268e+00& 8.04e+01 \\
nblk= 3 & {\bf SDPNAL+} & 4.91e-07 & 1.19e-09 & 6.86e-10 & -1.7921268e+00& 6.17e+01 \\
m=5051 & SDPT3 & 2.86e-10 & 0.00e+00 & 4.55e-08 & -1.7921263e+00& 8.78e+01 \\ 
\hline
\end{longtable}
\end{tiny}
\end{center}

\end{small}

\end{document}